\documentclass[12pt]{amsart}

\usepackage{geometry}
\usepackage{amsfonts, adjustbox, amssymb, amscd}
\numberwithin{equation}{section}

\usepackage{bm}
\usepackage{verbatim}
\usepackage{mathrsfs}
\usepackage{graphicx}
\usepackage{tikz-cd}
\usepackage{subcaption}
\usepackage{listings}
\usepackage{subfiles}
\usepackage[toc,page]{appendix}
\usepackage{mathtools}
\usepackage{comment}
\usepackage{enumerate}
\usepackage{enumitem}
\usepackage[all]{xy}
\usepackage{calligra}

\usepackage{graphicx}
\usepackage{appendix}
\usepackage{hyperref}
\hypersetup{
    colorlinks=true,
    citecolor=red,
    linkcolor=blue,
    filecolor=magenta,      
    urlcolor=red,
}
\newcommand{\bb}{\bm{b}}

\newcommand{\Mm}{{\bf{M}}}
\newcommand{\Nn}{{\bf{N}}}

\newcommand{\Pp}{{\bf{P}}}

\newcommand{\Qq}{\mathbb{Q}}

\newcommand{\Rr}{\mathbb{R}}

\newcommand{\ff}{\mathfrak{f}}

\newcommand{\vol}{\operatorname{vol}}

\newcommand{\Exc}{\operatorname{Exc}}

\newcommand{\rk}{\operatorname{rank}}

\newcommand{\ninv}{\operatorname{ninv}}
\newcommand{\inv}{\operatorname{inv}}

\newcommand{\tmld}{{\operatorname{tmld}}}

\DeclareMathOperator{\HHom}{\mathscr{H}\text{\kern -3pt {\calligra\large om}}\,}

\newcommand{\Supp}{\operatorname{Supp}}

\newcommand{\mult}{\operatorname{mult}}

\newcommand{\Aa}{{\mathfrak{A}}}

\newcommand{\Bb}{\mathcal{B}}
\newcommand{\Ff}{\mathcal{F}}

\newcommand{\Ii}{\Gamma}


\newcounter{parentnumber}

\newtheorem{thm}{Theorem}[section]

\newtheorem{conj}[thm]{Conjecture}
\newtheorem{cor}[thm]{Corollary}
\newtheorem{lem}[thm]{Lemma}
\newtheorem{prop}[thm]{Proposition}

\theoremstyle{definition}
\newtheorem{defn}[thm]{Definition}

\theoremstyle{definition}
\newtheorem{rem}[thm]{Remark}

\newtheorem{defthm}[thm]{Definition-Theorem}

\newtheorem{nota}[thm]{Notation}

\theoremstyle{definition}

\begin{document}

\title{Non-algebraicity of non-abundant foliations and abundance for adjoint foliated structures}

\author{Jihao Liu}
\author{Zheng Xu}

\date{\today}

\begin{abstract}
Assuming the abundance conjecture in dimension $d$, we establish a non-algebraicity criterion of foliations: any log canonical foliation of rank $\le d$ with $\nu\neq\kappa$ is not algebraically integrable, answering question of Ambro--Cascini--Shokurov--Spicer. Under the same hypothesis, we prove abundance for klt algebraically integrable adjoint foliated structures of dimension $\le d$ and show the existence of good minimal models or Mori fiber spaces. In particular, when $d=3$, all these results hold unconditionally. 

Using similar arguments, we solve a problem proposed by Lu and Wu on abundance of surface adjoint foliated structures that are not necessarily algebraically integrable.
\end{abstract}

\address{Department of Mathematics, Peking University, No. 5 Yiheyuan Road, Haidian District, Beijing 100871, China}
\email{liujihao@math.pku.edu.cn}

\address{Department of Mathematics, Southern University of Science and Technology (SUSTech), 1088 Xueyuan Rd, Nanshan District, Shenzhen, Guangdong, China, 518055}
\email{xuz6@sustech.edu.cn}

\maketitle

\pagestyle{myheadings}\markboth{\hfill Jihao Liu and Zheng Xu \hfill}{\hfill Non-algebraicity of non-abundant foliations and abundance for adjoint foliated structures\hfill}

\tableofcontents

\section{Introduction}\label{sec:Introduction}

We work over the field of complex numbers $\mathbb{C}$.

This paper studies two central problems in the study of foliations on complex manifolds: the \emph{algebraicity} of the foliation and the \emph{birational classification} of foliations.

\medskip

\noindent\textbf{Algebraicity of foliations.} The algebraicity problem seeks criteria for the algebraic integrability of a foliation, i.e., when the general leaves of foliations are algebraic. For foliated surfaces, this question is also known as the ``Poincaré problem'' and its study can be traced back to the 19th century, when Darboux (\cite{Dar78a,Dar78b}), Poincaré (\cite{Poi85,Poi91}), Painlevé (\cite{Pai92}), etc. studied the orbits of vector fields on $\mathbb P^2$. A famous theorem of Darboux–Jouanolou–Ghys \cite{Dar78a,Jou78,Ghy00} indicates that a surface foliation $\Ff$ is algebraically integrable as long as it admits sufficiently many invariant algebraic curves. Moreover, it is known that a surface foliation $\Ff$ with non-pseudo-effective canonical divisor $K_{\Ff}$ is uniruled, hence algebraically integrable \cite{Miy87}. A significant advancement on the Poincaré problem was recently obtained by Lü and Tan \cite{LT24}, where new birational invariants on foliated surfaces, motivated by the formulas for the modular invariants of families of algebraic curves (cf. \cite{Tan96,Tan10}), were introduced, and several criteria for the algebraicity of surface foliations were established. 

However, the algebraicity problem is substantially more difficult in higher dimensions. After \cite{Miy87}, the uniruled property for foliations with non-pseudo-effective canonical divisor was generalized to higher dimensions \cite{Bos01,BM16,CP19}, and very recently to the Kähler setting \cite{Ou25a,CP25}. There are also significant advancements on the algebraicity of foliations with numerically trivial or anti-nef canonical bundles in recent years \cite{LPT18,Dru18,HP19,Dru21,DO22,Ou25b}. Yet, we know very little about the algebraicity of foliations with “positive’’ (but not trivial) canonical classes.

\medskip

\noindent\textbf{Birational classification of foliations.} The birational classification of foliations aims to classify foliations up to birational equivalence by considering the behavior of their birational invariants. The classification is in the same spirit as the Enriques classification of algebraic surfaces (going back to the Italian school) and the minimal model program: we first run a foliated minimal model program for the foliation, obtain a minimal model (or canonical model), and classify the foliations according to the Kodaira dimension $\kappa(\Ff)$ and the numerical dimension $\nu(\Ff)$. For foliated surfaces, thanks to Brunella \cite{Bru15} and McQuillan \cite{McQ08}, the minimal model program and such a classification have been completed, and the canonical models (in the sense of \cite[Definition III.3.1]{McQ08}) of all non-general type canonical foliations on surfaces are explicitly described. In more recent works, the existence of the minimal model program has been established for threefolds (cf. \cite{CS25c,CS21,SS22}), algebraically integrable foliations (cf. \cite{ACSS21,CHLX23,LMX25,CS25a}), and partially for rank one foliations in arbitrary dimensions (cf. \cite{McQ24,CS25b,CS26}). See also \cite{CM24,Li25}.

Despite the aforementioned progress, one major issue for the classification of foliations is the failure of abundance, even for surfaces. More precisely:
\begin{enumerate}
    \item (cf. \cite[Theorem IV.2.2]{McQ08}) There exists a canonical foliation $\Ff$ on a canonical projective surface $X$ such that $K_{\Ff}$ is big and nef but $K_{\Ff}$ is not semi-ample. In particular, the canonical ring 
$$R(X,K_{\Ff}):=\oplus_{m=0}^{+\infty}H^0(X,\mathcal{O}_X(mK_{\Ff}))$$
    is not finitely generated.
    \item (cf. \cite[Theorem IV.5.11]{McQ08}) There exists a canonical foliation $\Ff$ on a canonical projective surface $X$ such that $K_{\Ff}$ is nef but
    $$1=\nu(\Ff)\not=\kappa(\Ff)=-\infty.$$
\end{enumerate}

This is very different from the case of the birational geometry of varieties, as the finite generation and base-point-freeness theorem are both known (cf. \cite[Theorem 1.2]{BCHM10}, \cite[Theorem 3-1-1]{KMM87}), while the abundance conjecture (Conjecture \ref{conj: abundance}) predicts that $\nu(X)=\kappa(X)$ and has been proved in dimension $\leq 3$ (cf. \cite{Kaw92,KMM94}). 

\medskip

\noindent\textbf{Main Theorems.} The first main theorem of this paper is the following non-algebraicity criterion for foliations that are not abundant, i.e., $\nu(\Ff)\neq\kappa(\Ff)$.

\begin{thm}\label{thm: main nonalg I}
     Assume the abundance conjecture (Conjecture \ref{conj: abundance}) in dimension $d$. Let $\Ff$ be a log canonical foliation of rank $\leq d$ such that $\nu(\Ff)\neq\kappa(\Ff)$. 
     
     Then $\Ff$ is not algebraically integrable.
\end{thm}
We refer the reader to Theorem \ref{thm: log abundance aif} for a more general result.

Theorem \ref{thm: main nonalg I} provides a positive answer to a question of Ambro, Cascini, Shokurov, and Spicer \cite[Subsection 5.4]{ACSS21} assuming the abundance conjecture. It also provides a new (conditional) non-algebraicity criterion: Note that, by the classification of surface foliations \cite{Bru15,McQ08}, non-abundant canonical foliations form a non-empty (although possibly small) class. In particular, we can show that surface canonical foliations $\Ff$ with $\nu(\Ff)=1$ and $\kappa(\Ff)=-\infty$ are not algebraically integrable. This is already known via the classification \cite[Theorem IV.5.11(3)]{McQ08}, yet by applying Theorem \ref{thm: main nonalg I} we obtain a proof without classification. Moreover, since the abundance conjecture holds in dimension $\leq 3$ \cite[1.1 Theorem]{KMM94}, we immediately obtain the following:

\begin{cor}\label{cor: rank 3 abudnance}
 Let $\Ff$ be a log canonical foliation of rank $\leq 3$ such that $\nu(\Ff)\neq\kappa(\Ff)$. Then $\Ff$ is not algebraically integrable.   
\end{cor}

Our next theorem concerns the existence of good minimal models of foliations. For smooth projective varieties $X$, the abundance conjecture, which predicts that $\nu(X)=\kappa(X)$, implies the existence of good minimal model or Mori fiber spaces of $X$ (cf. \cite[Theorem 4.3]{GL13}, \cite[Theorem 4.4]{Lai11}). However, this implication fails for foliations, even for algebraically integrable ones (cf. \cite[Theorem IV.2.2]{McQ08}, \cite[Examples 5.4, 5.5]{ACSS21}). This is a big problem as we can no longer consider the ample models foliations and hence cannot deduce a nice boundedness and moduli theory (cf. \cite[Example 4.1]{Pas24}).

One very recent approach to tackle this issue is to consider \emph{adjoint foliated structures}, i.e., instead of considering $K_{\Ff}$, we consider the structure $(X,\Ff,t)$ associated with the canonical divisor $K_t:=tK_{\Ff}+(1-t)K_X$. A lot of works have recently been done towards the birational geometry of adjoint foliated structures \cite{PS19,SS23,Cas+24,LW24,Cas+25a,CC25,LWX26}, and now there is substantial evidence that the birational geometry of $K_t$ (with $t<1$) behaves much better than the birational geometry of $K_{\Ff}$ and is closer to that of $K_X$. The most important piece of evidence is the existence of good minimal models, and hence the finite generation of the canonical ring, for adjoint foliated structures of general type in arbitrary dimensions \cite[Theorem A]{Cas+25a}.

In this paper, we show that abundance holds for algebraically integrable adjoint foliated structures and that good minimal models and Mori fiber spaces always exist, provided that the usual abundance conjecture holds: 
\begin{thm}\label{thm: main abundance aiafs}
     Assume the abundance conjecture (Conjecture \ref{conj: abundance}) in dimension $d$. Let $X$ be a smooth projective variety of dimension $\leq d$ and let $\Ff$ be a log canonical algebraically integrable foliation on $X$. We let $K_t:=tK_{\Ff}+(1-t)K_X$ for any $t\in [0,1]$. Then:
     \begin{enumerate}
         \item For any $t\in [0,1]$, $\nu(K_t)=\kappa_{\iota}(K_t)$. Here $\kappa_{\iota}$ stands for the invariant Iitaka dimension (Definition \ref{defn: iitaka dimension}) and we have $\kappa_{\iota}(K_t)=\kappa(K_t)$ when $t\in\mathbb Q$.
         \item For any $t\in [0,1)$, we may run a $K_t$-MMP which terminates with a good minimal model or a Mori fiber space of $(X,\Ff,t)$.
     \end{enumerate}
\end{thm}
We refer the reader to Theorem \ref{thm: log abundance aiafs assume abundance} for a more general result.

It is interesting to ask whether abundance is expected to hold for adjoint foliated structures that are not necessarily algebraically integrable. Based on the proof of Theorem \ref{thm: main abundance aiafs}, adjoint foliated structures $(X,\Ff,t)$ with abundant $K_{\Ff}$ can be treated in a similar way modulo the necessary MMPs, and the main difficulty comes from the non-algebraically integrable ones. For surfaces, however, thanks to the classification of non-abundant foliations \cite[Theorem IV.5.11]{McQ08}, we are able to deduce abundance in full generality:

\begin{thm}\label{thm: abundance afs surface}
Let $X$ be a smooth projective surface, $\Ff$ a canonical foliation on $X$, $t\in [0,1)$ a real number, and $K_t:=tK_{\Ff}+(1-t)K_X$. Then $\nu(K_t)=\kappa_{\iota}(K_t)$, and if $K_t$ is pseudo-effective, then $(X,\Ff,t)$ has a good minimal model.
\end{thm}
We remark that \cite[Problem 1.9]{LW24} asks whether there exists a relatively minimal foliation $\Ff$ on a smooth projective surface $X$ such that $\kappa(K_X+K_{\Ff})\not=\nu(K_X+K_{\Ff})$. Since any relatively minimal foliation is always canonical, Theorem \ref{thm: abundance afs surface} provides a negative answer to \cite[Problem 1.9]{LW24} by considering the special case of $t=\frac{1}{2}$.

\medskip

\noindent\textit{Sketch of the proof.} We first prove Theorem \ref{thm: main nonalg I}. For simplicity, assume that $\Ff$ is canonical. By \cite{ACSS21,CHLX23} we may reduce to the case where $\Ff$ is induced by an equidimensional contraction $f\colon X\rightarrow Z$, $X$ is BP stable$/Z$, and $K_{\Ff}$ is semi-ample$/Z$ and nef. After taking a finite cover, we may further reduce to the case where $f$ is a locally stable family, hence $K_{\Ff}=K_{X/Z}$. Let $g\colon X\rightarrow Y$ be the ample model$/Z$ of $K_{\Ff}$. Then
$$
K_{\Ff}\sim_{\mathbb Q}g^*(K_{\Ff_Y}+B_Y+\Mm_Y)
$$
where $(Y,\Ff_Y,B_Y,\Mm)$ is a generalized foliated quadruple with $K_{\Ff_Y}=K_{Y/Z}$ and $h\colon (Y,B_Y,\Mm)\rightarrow Z$ is the induced stable family of generalized pairs. 

If $\Mm=\bm{0}$, then by \cite[Corollary 6.20]{KP17} the $\mathbb R$-divisor $K_{Y/Z}+B_Y$ is crepant to $K_{T/Z'}+B_T$, where $(T,B_T)\rightarrow Z'$ is a stable family of maximal variation with projective base $Z'$ and lc general fiber. By \cite[Proposition 2.15]{PX17} we have that $K_{T/Z'}+B_T$ is big; hence $K_{Y/Z}+B_Y$ is abundant, and so is $K_{\Ff}$. For the general case, following \cite[Proof of Theorem 4.6]{FS25}, we apply \cite[Theorem 1.5]{BFMT25} to choose a finite open cover $\{U_i\}$ of $Z$, let $V_i$ be the inverse image of $U_i$ in $Y$, and choose $L_i\in |\Mm_Y|_{\mathbb Q}$ so that $(Y,B_Y+L_i)|_{V_i}\rightarrow U_i$ is a stable family of lc pairs.

For any prime divisor $P$ over $Y$, choose $V=V_i$ containing the generic point of its center. After finite base change, we obtain a stable family $(Y',D')\rightarrow Z'$ with projective base and normal general fiber, and a smooth projective variety $W$ with morphisms $r:W\rightarrow Y$ and $q:W\rightarrow Y'$, where $r$ is generically finite and surjective and $q$ is birational. For $H:=q^*(K_{Y'/Z'}+D')$, there is an effective $\mathbb R$-divisor $E$ on $W$ such that
\[
r^*(K_{Y/Z}+B_Y+\Mm_Y)\sim_{\mathbb R}H+E,
\qquad E\geq0,\qquad E|_{r^{-1}(V)}=0,
\]
where $H$ is nef and abundant by Lemma~\ref{lem: abudant stable family}.
Put $D:=K_{Y/Z}+B_Y+\Mm_Y$, and let $\mu:Y_0\rightarrow Y$ be a projective birational morphism with $Y_0$ normal and $P$ a prime divisor on $Y_0$. The calculation in Proposition~\ref{prop: nef relative sa imply abundant} gives, for every $\epsilon>0$, an effective $\mathbb R$-divisor $N_{P,\epsilon}\sim_{\mathbb R}\mu^*D$ with $\mult_PN_{P,\epsilon}<\epsilon$. Since $g^*D\sim_{\mathbb Q}K_{\Ff}$ is nef and $g$ is surjective, $D$ is nef. Lemma~\ref{lem: equivalent definition of abundance}(3) therefore implies that $D$ is abundant, and hence so is $K_{\Ff}$.

To deduce Theorems \ref{thm: main abundance aiafs} and \ref{thm: abundance afs surface}, we adapt the idea of \cite[Proof of Theorem 5.2]{Cas+25a}, with minor modifications. We may assume that $K_t$ is pseudo-effective. Let $\Aa_s:=(X,\Ff,s)$ for any $s\in [0,1]$. Consider the sets
$$
\mathcal{P}:=\{s\in [t,1]\mid \Aa_s\text{ has a }\mathbb Q\text{-factorial good minimal model}\}
$$
and
$$
\mathcal{Q}:=\{s\in [t,1]\mid \Aa_s\text{ has a }\mathbb Q\text{-factorial minimal model and }K_{\Aa_s}\text{ is abundant}\}.
$$
We prove the following:
\begin{itemize}
    \item[(i)] If $\mu\in\mathcal{Q}$ and $\mu>t$, then there exists $\epsilon>0$ such that $[\mu-\epsilon,\mu)\subset\mathcal{P}$.
    \item[(ii)] If $\mu<1$ and $\mu+\epsilon\in\mathcal{P}$ for any $0<\epsilon\ll 1$, then $\mu\in\mathcal{P}$.
\end{itemize}
The key input for both (i) and (ii) is the following observation: if $K_{\Aa_s}$ is nef and abundant for some $s$, then for any $s'<s$,
$$
\frac{s}{s-s'}K_{\Aa_{s'}}\sim_{\mathbb R}K_{\Aa_0}+\frac{s'}{s-s'}K_{\Aa_s}
$$
admits a klt generalized pair structure whose nef part is nef and abundant, and hence underlies a klt pair structure (cf. Lemma \ref{lem: klt with nef and abundant boundary}). In this situation, the usual abundance conjecture implies the abundance of $K_{\Aa_{s'}}$. See Theorem \ref{thm: log smooth existence of good minimal model} for details.

Combining (i) and (ii) with the inclusion $\mathcal{P}\subset\mathcal{Q}$, we deduce that if $\mathcal{Q}\cap (t,1]\not=\emptyset$, then $t\in\mathcal{P}$, and we are done by using the fact that the existence of good minimal models implies the termination of the MMP with scaling (cf. \cite[Theorem 4.5]{Cas+25b}). 

Theorem \ref{thm: main abundance aiafs} now follows immediately from Theorem \ref{thm: main nonalg I}, since $1\in\mathcal{Q}$. The proof of Theorem \ref{thm: abundance afs surface} is similar. If $1\in\mathcal{Q}$ then we are done. If $1\not\in\mathcal{Q}$, write $K_{\Ff}=P+N$ for its Zariski decomposition on $X$. By \cite[Lemma IV.3.1, Proposition IV.5.1, and Lemma IV.5.2]{McQ08}, we have $P^2=0$, $\kappa(K_{\Ff})=-\infty$, and $q(X)=0$. Riemann--Roch gives $P\cdot K_X>0$, so $sP+(1-s)K_X$, and hence $K_{\Aa_s}$, is big for $s<1$ sufficiently close to $1$. Thus such an $s$ lies in $\mathcal Q$, as required. Of course, several preparatory results are needed for non-algebraically integrable adjoint foliated structures on surfaces (cf. Subsection \ref{subsec: surface afs}).

\medskip

\noindent\textit{Structure of the paper.} In Section \ref{sec: preliminaries} we recall preliminary material and recall useful results we need later. In Section \ref{sec: abu f} we prove Theorem \ref{thm: main nonalg I}. In Section \ref{sec: abu aiafs} we prove Theorem \ref{thm: main abundance aiafs}. Log versions of Theorems \ref{thm: main nonalg I} and \ref{thm: main abundance aiafs} are also established in the corresponding sections. In Section \ref{sec: abu surface}, we first develop auxiliary results on surface adjoint foliated structures, and then prove Theorem \ref{thm: abundance afs surface}.

\section{Preliminaries}\label{sec: preliminaries}

We will adopt the standard notation and definitions on MMP in \cite{Sho92,KM98,BCHM10} and use them freely. For adjoint foliated structures, generalized foliated quadruples, foliated triples, foliations, and generalized pairs, we adopt the notation and definitions in \cite{Cas+24,Cas+25a}  which generally align with \cite{LLM23,CHLX23} (for generalized foliated quadruples), \cite{CS25c,ACSS21,CS21} (for foliations and foliated triples), and \cite{BZ16,HL23} (for generalized pairs and $\bb$-divisors), possibly with minor differences.  We refer the reader to \cite{Kol23} for basic notations and properties of moduli theory but we shall recall some important concepts here.

\subsection{Basic notations}

\begin{defn}
    Let $\Ii\subset [0,1]$ be a set. We define
    $$\Ii_+:=\{0\}\cup\left(\left\{\sum_{i=1}^n\gamma_i \mid n\in\mathbb N^+,\gamma_1,\dots,\gamma_n\in\Ii\right\}\cap [0,1]\right).$$
\end{defn}

\begin{defn}
    A \emph{contraction} is a projective morphism $f: X\rightarrow Y$ such that $f_*\mathcal{O}_X=\mathcal{O}_Y$. A \emph{finite cover} is a finite surjective morphism. 
\end{defn}

\begin{nota}
    Let $f: X\dashrightarrow X'$ be a birational map between normal schemes. We denote by $\Exc(f)$ the reduced divisor supported on the codimension one part of the exceptional locus of $f$.
\end{nota}

\begin{defn}
Let $m$ be a positive integer and $\bm{v}\in\mathbb R^m$. The \emph{rational envelope} of $\bm{v}$ is the minimal rational affine subspace of $\mathbb R^m$ which contains $\bm{v}$. For example, if $m=2$ and $\bm{v}=\left(\frac{\sqrt{2}}{2},1-\frac{\sqrt{2}}{2}\right)$, then the rational envelope of $\bm{v}$ is $(x_1+x_2=1)\subset\mathbb R^2_{x_1x_2}$.
\end{defn}

\begin{defn}
    Let $X\rightarrow U$ be a projective morphism from a normal quasi-projective variety to a variety.  Let $D$ be an $\Rr$-Cartier $\Rr$-divisor on $X$ and $\phi: X\dashrightarrow X'$ a birational map$/U$. We say that $\phi$ is \emph{$D$-negative} (resp. \emph{$D$-trivial}) if the following conditions hold:
    \begin{enumerate}
    \item $\phi$ does not extract any divisor.
    \item $D':=\phi_\ast D$ is $\Rr$-Cartier.
    \item There exists a resolution of indeterminacy $p: W\rightarrow X$ and $q: W\rightarrow X'$, such that
    $$p^\ast D=q^\ast D'+F$$
    where $F\geq 0$ and $\Supp p_\ast F=\Exc(\phi)$ (resp. $F=0$).
    \end{enumerate}
\end{defn}

\subsection{Foliation and adjoint foliated structures}

\begin{defn}[Foliations, {cf. \cite{ACSS21,CS21}}]\label{defn: foliation}
Let $X$ be a normal variety. A \emph{foliation} on $X$ is a coherent subsheaf $T_{\Ff}\subset T_X$ such that
\begin{enumerate}
    \item $T_{\Ff}$ is saturated in $T_X$, i.e. $T_X/T_{\Ff}$ is torsion free, and
    \item $T_{\Ff}$ is closed under the Lie bracket.
\end{enumerate}
The \emph{canonical divisor} of $\Ff$ is a divisor $K_\Ff$ such that $\mathcal{O}_X(-K_{\mathcal{F}})\cong\mathrm{det}(T_\Ff)$. If $T_\Ff=0$, then we say that $\Ff$ is a \emph{foliation by points}.

Given any dominant map 
$h: Y\dashrightarrow X$ and a foliation $\mathcal F$ on $X$, we denote by $h^{-1}\Ff$ the \emph{pullback} of $\Ff$ on $Y$ as constructed in \cite[3.2]{Dru21}. Given any birational map $g: X\dashrightarrow X'$, we denote by $g_\ast \Ff:=(g^{-1})^{-1}\Ff$ the \emph{pushforward} of $\Ff$ on $X'$. We say that $\Ff$ is an \emph{algebraically integrable foliation} if there exists a dominant map $f: X\dashrightarrow Z$ such that $\Ff=f^{-1}\Ff_Z$, where $\Ff_Z$ is the foliation by points on $Z$, and we say that $\Ff$ is \emph{induced by} $f$.

A subvariety $S\subset X$ is called \emph{$\Ff$-invariant} if for any open subset $U\subset X$ and any section $\partial\in H^0(U,\Ff)$, we have $\partial(\mathcal{I}_{S\cap U})\subset \mathcal{I}_{S\cap U}$,  where $\mathcal{I}_{S\cap U}$ denotes the ideal sheaf of $S\cap U$ in $U$.  
For any prime divisor $P$ on $X$, we define $\epsilon_{\Ff}(P):=1$ if $P$ is not $\Ff$-invariant and $\epsilon_{\Ff}(P):=0$ if $P$ is $\Ff$-invariant. For any prime divisor $E$ over $X$, we define $\epsilon_{\Ff}(E):=\epsilon_{\Ff_Y}(E)$ where $h: Y\dashrightarrow X$ is a birational map such that $E$ is on $Y$ and $\Ff_Y:=h^{-1}\Ff$.

Suppose that the foliation structure $\Ff$ on $X$ is clear in the context. Then, given an $\mathbb R$-divisor $D =\sum a_iD_i$ where each $D_i$ is a prime divisor,
we denote by $D^{\ninv} \coloneqq \sum \epsilon_{\Ff}(D_i)a_iD_i$ and $D^{\inv} \coloneqq D-D^{\ninv}$.
\end{defn}

\begin{defn}\label{defn: afs}
An \emph{adjoint foliated structure} $\Aa/U:=(X,\Ff,B,\Mm,t)/U$ is the datum of a normal quasi-projective variety $X$ and a projective morphism $X\rightarrow U$, a foliation $\Ff$ on $X$, an $\Rr$-divisor $B\geq 0$ on $X$, a nef$/U$ $\bb$-divisor $\Mm$, and a real number $t\in [0,1]$ such that $K_{\Aa}:=tK_{\Ff}+(1-t)K_X+B+\Mm_X$ is $\Rr$-Cartier. We may simply say that ``$\Aa/U$ is an adjoint foliated structure" without mentioning $X,\Ff,B,\Mm,t$ at all. $X,t$ are called the \emph{ambient variety} and \emph{parameter} of $\Aa$ respectively. 

We say that $\Aa/U$ is \emph{of general type} if $K_{\Aa}$ is big$/U$. We say that $\Aa$ is \emph{algebraically integrable} if $\Ff$ is algebraically integrable.

For any $\Rr$-divisor $D$ on $X$ and nef$/U$ $\bb$-divisor $\Nn$ on $X$ such that $D+\Nn_X$ is $\Rr$-Cartier, we denote by $(\Aa,D,\Nn):=(X,\Ff,B+D,\Mm+\Nn,t)$. If $D=0$ then we may drop $D$, and if $\Nn=\bm{0}$ then we may drop $\Nn$.

When $t=0$ or $\Ff=T_X$, we call $(X,B,\Mm)/U$ a \emph{generalized pair}, and in addition, if $\Mm=\bm{0}$, then we call $(X,B)/U$ a pair. If $B=0$, or if $\Mm=\bm{0}$, or if $U$ is not important, then we may drop $B,\Mm,U$ respectively. If $U=\{pt\}$ then we also drop $U$ and say that $(X,\Ff,B,\Mm,t)$ is \emph{projective}.

For any birational map$/U$ $\phi: X\dashrightarrow X'$, we define $\phi_*\Aa:=(X',\phi_*\Ff,\phi_*B,\Mm,t)$ and say that $\phi_*\Aa$ is the \emph{image} of $\Aa$ on $X'$. 
For any projective birational morphism $h: X'\rightarrow X$, we define 
$$h^*\Aa:=(X',\Ff',B',\Mm,t)$$
where $\Ff':=h^{-1}\Ff$ and $B'$ is the unique $\Rr$-divisor such that $K_{h^*\Aa}=h^*K_{\Aa}$. For any prime divisor $E$ on $X'$, we denote by
$$a(E,\Aa):=-\mult_EB'$$
the \emph{discrepancy} of $E$ with respect to $\Aa$. The \emph{total minimal log discrepancy} of $\Aa$ is
$$\tmld(\Aa):=\inf\{a(E,\Aa)+t\epsilon_{\Ff}(E)+(1-t)\mid E\text{ is over }X\}.$$
For any non-negative real number $\epsilon$, we say that $\Aa$ is \emph{$\epsilon$-lc} (resp. \emph{$\epsilon$-klt}) if $\tmld(\Aa)\geq\epsilon$ (resp. $>\epsilon$). We say that $\Aa$ is \emph{lc} (resp. \emph{klt}) if $\Aa$ is $0$-lc (resp. $0$-klt). We say that $\Aa$ is \emph{canonical} if $a(E,\Aa)\geq 0$ for any prime divisor $E$ that is exceptional$/X$. 
\end{defn}

\begin{defn}[Potentially klt]\label{defn: potentially klt}
Let $X$ be a normal quasi-projective variety. We say that $X$ is \emph{potentially klt} if $(X,\Delta)$ is klt for some $\Rr$-divisor $\Delta\geq 0$. 
\end{defn}

\begin{defn}[{\cite[Definition 3.11]{Cas+25a}}]\label{defn: foliated log smooth}
Let $\Aa/U:=(X,\Ff,B,\Mm,t)/U$ be an algebraically integrable adjoint foliated structure. We say that $\Aa$ is \emph{foliated log smooth} if there exists a contraction $f: X\rightarrow Z$ satisfying the following.
\begin{enumerate}
  \item $X$ has at most quotient singularities.
  \item $\Ff$ is induced by $f$.
  \item $(X,\Sigma_X)$ is toroidal for some reduced divisor $\Sigma_X$ such that $\Supp B\subset\Sigma_X$.  In particular, $(X,\Supp B)$ is toroidal, and $X$ is $\Qq$-factorial klt.
  \item There exists a log smooth pair $(Z,\Sigma_Z)$ such that $$f: (X,\Sigma_X)\rightarrow (Z,\Sigma_Z)$$ is an equidimensional toroidal contraction.
  \item $\Mm$ descends to $X$.
\end{enumerate}
We say that $f: (X,\Sigma_X)\rightarrow (Z,\Sigma_Z)$ is \emph{associated with} $\Aa$.
\end{defn}

\begin{defn}[Foliated log resolutions]\label{defn: log resolution}
Let $\Aa/U:=(X,\Ff,B,\Mm,t)/U$ be an algebraically integrable adjoint foliated structure. A \emph{foliated log resolution} of $\Aa$ is a birational morphism $h: X'\rightarrow X$ such that $(h^{-1}_*\Aa,\Exc(h))$ is foliated log smooth. By \cite[Lemma 6.2.4]{CHLX23}, a foliated log resolution for $\Aa$ always exists.
\end{defn}

\begin{defn}[Property $(\ast )$ foliations, {\cite[Definition 3.8]{ACSS21}, \cite[Definition 7.2.2]{CHLX23}}]\label{defn: foliation property *}
Let $\Aa/U:=(X,\Ff,B,\Mm)/U$ be a generalized foliated quadruple, $G\geq 0$ be a reduced divisor on $X$, and $f: X\rightarrow Z$ a contraction. We say that $(\Aa;G)/Z$ satisfies \emph{Property $(\ast )$} if the following conditions hold.
\begin{enumerate}
\item $\Ff$ is induced by $f$ and $G$ is an $\Ff$-invariant divisor.
\item $f(G)$ is of pure codimension $1$, $(Z,f(G))$ is log smooth, and $G=f^{-1}(f(G))$.
\item For any closed point $z\in Z$ and any reduced divisor  $\Sigma\ge f(G)$ on $Z$ such that  $(Z,\Sigma)$ is log smooth near $z$, $(X,B+G+f^\ast (\Sigma-f(G)),\Mm)$ is lc over a neighborhood of $z$.
\end{enumerate}
We say that $f$, $Z$, and $G$ are \emph{associated} with $\Aa$. 
\end{defn}

\begin{defn}[ACSS, {cf. \cite[Definitions 5.4.2, 7.2.2, 7.2.3]{CHLX23}}]\label{defn: ACSS f-triple}
Let $\Aa/U:=(X,\Ff,B,\Mm)/U$ be an lc generalized foliated quadruple, $G\geq 0$ a reduced divisor on $X$, and $f: X\rightarrow Z$ a contraction. We say that $(\Aa;G)/Z$ is \emph{ACSS} if the following conditions hold:
\begin{enumerate}    
\item $(\Aa;G)/Z$ satisfies Property $(\ast )$.
\item $f$ is equidimensional.
\item There exists an $\Rr$-Cartier $\Rr$-divisor $D\geq 0$ on $X$ and a nef$/X$ $\bb$-divisor $\Nn$ on $X$, such that  $\Supp\{B\}\subset\Supp D$, $\Nn-\alpha\Mm$ is nef$/X$ for some $\alpha>1$, and for any reduced divisor $\Sigma\geq f(G)$ such that $(Z,\Sigma)$ is log smooth, $$(X,B+D+G+f^\ast (\Sigma-f(G)),\Nn)$$ 
      is qdlt (cf. \cite[Definition 7.1.1]{CHLX23}, \cite[Definition 35]{FKX17}).
\item For any lc center of $\Aa$ with generic point $\eta$, over a neighborhood of $\eta$,
    \begin{enumerate}
      \item $\eta$ is the generic point of an lc center of $(X,\Ff,\lfloor B\rfloor)$, and
       \item $f: (X,B+G)\rightarrow (Z,f(G))$ is a toroidal morphism.
    \end{enumerate}
\end{enumerate}
If $(\Aa;G)/Z$ is ACSS, then we say that $\Aa/Z$ and $\Aa$ are ACSS.
\end{defn}

\begin{defn}\label{def: relative ample model}
    Let $X',X,Z$ be normal quasi-projective varieties and $h: X'\rightarrow X$, $f: X'\rightarrow Z$ contractions. The \emph{core model} of $(h,f)$ associated with $(\bar h,\bar f)$ is the unique normal quasi-projective variety $\bar X$ up to isomorphism with two contractions $\bar h: \bar X\rightarrow X$ and $\bar f: \bar X\rightarrow Z$ satisfying the following.
    \begin{enumerate}
        \item  For any ample $\Rr$-divisor $A$ on $X$, $\bar h^*A$ is ample$/Z$.
        \item  There exists a contraction $g: X'\rightarrow \bar X$ such that $\bar h\circ g=h$ and $\bar f\circ g=f$.
    \end{enumerate}
The variety $\bar X$ is called the \emph{core model} of $(h,f)$ associated with $(\bar h,\bar f)$. Existence of core model is guaranteed by \cite[Definition-Lemma 3.1]{LMX25}.
\end{defn}

\begin{defn}[{\cite[Subsection 3.2]{LMX25}}]\label{deflem: simple model}
    Let $\Aa:=(X,\Ff,B,\Mm)$ and $\Aa':=(X',\Ff',B',\Mm)$ be two algebraically integrable generalized foliated quadruples and $h: X'\rightarrow X$ a birational morphism. Let $f: X'\rightarrow Z$ be a contraction and $G$ a reduced divisor on $X'$. We say that $h: (\Aa';G)/Z\rightarrow \Aa$ is a \emph{simple modification} if the following conditions hold.
    \begin{itemize}
        \item $\Aa'=(h^{-1}_*\Aa,\Exc(h)^{\ninv})$ and $\Aa'$ is lc.
        \item  $a(E,\Aa)\leq-\epsilon_{\Ff}(E)$ for any $h$-exceptional prime divisor $E$. 
        \item $K_{\Aa'}\sim_ZK_{X'}+B'+\Mm_{X'}+G.$
        \item $(\Aa';G)/Z$ satisfies Property $(*)$.
    \end{itemize}
    We say that $h: (X',\Ff',B',G)/Z\rightarrow\Aa$ is 
    \begin{enumerate}
        \item an \emph{ACSS modification} if it is a simple modification and  $(\Aa';G)/Z$ is ACSS, and we say that $\Aa'$ is an \emph{ACSS model} of $\Aa$. 
        \item a \emph{core modification} if it is a simple modification and $h^*A$ is ample$/Z$ for any ample $\Rr$-divisor $A$ on $X$, and we say that $\Aa'$ is a \emph{core model} of $\Aa$, and
        \item \emph{$\Qq$-factorial} if $X'$ is $\Qq$-factorial.
    \end{enumerate}
\end{defn}

\subsection{Models in the minimal model program}

\begin{defn}[Log birational models]\label{defn: log birational model}
Let $\Aa/U$ be an adjoint foliated structure with ambient variety $X$, $\phi: X\dashrightarrow X'$ a birational map$/U$, and $E:=\Exc(\phi^{-1})$ the reduced $\phi^{-1}$-exceptional divisor. Assume that $a(D,\Aa)\leq-t\epsilon_{\Ff}(D)-(1-t)$ for any component $D$ of $E$. We let
  $$\Aa':=\left(\phi_*\Aa,-\sum_Da(D,\Aa)D\right)$$
where the sum runs through all irreducible components of $E$. If $K_{\Aa'}$ is $\mathbb R$-Cartier then we say that $\Aa'/U$ is a \emph{log birational model} of $\Aa/U$.
\end{defn}

\begin{defn}[Minimal models]\label{defn: minimal model}
    Let $\Aa/U$ be an adjoint foliated structure with ambient variety $X$ and $\Aa'/U$ a log birational model of $\Aa/U$ with ambient variety $X'$ and associated  birational map $\phi: X\dashrightarrow X'$, such that $K_{\Aa'}$ is nef$/U$. 

    Let $t$ be the parameter of $\Aa$.
    \begin{enumerate}
        \item We say that $\Aa'/U$ is a \emph{bs-weak lc model} or \emph{weak lc model in the sense of Birkar-Shokurov} of $\Aa/U$, if for any prime divisor $D$ on $X$ which is exceptional over $X'$, $$a(D,\Aa)\leq a(D,\Aa').$$
        We also say that $\phi$ is a bs-weak lc model of $\Aa/U$.
        \item We say that $\Aa'/U$ is a \emph{bs-minimal model} or \emph{minimal model in the sense of Birkar-Shokurov} of $\Aa/U$, if for any prime divisor $D$ on $X$ which is exceptional over $X'$, $$a(D,\Aa)<a(D,\Aa').$$
        We also say that $\phi$ is a bs-minimal model of $\Aa/U$.
        \item We say that $\Aa'/U$ is a \emph{bs-semi-ample model} or \emph{semi-ample model in the sense of Birkar-Shokurov} of $\Aa/U$ if it is a bs-weak lc model of $\Aa/U$ and $K_{\Aa'}$ is semi-ample$/U$.      We also say that $\phi$ is a bs-semi-ample model of $\Aa/U$.
        \item We say that $\Aa'/U$ is a \emph{bs-good minimal model} or \emph{good minimal model in the sense of Birkar-Shokurov} of $\Aa/U$ if it is a bs-minimal model of $\Aa/U$ and $K_{\Aa'}$ is semi-ample$/U$. We also say that $\phi$ is a bs-good minimal model of $\Aa/U$.
        \end{enumerate}
If, in addition, the induced birational map $X\dashrightarrow X'$ does not extract any divisor, then we remove the initial ``bs-" or the phrase ``in the sense of Birkar-Shokurov" in the previous definitions. 
\end{defn}

\begin{defn}[Mori fiber space]
Let $\Aa/U$ be an adjoint foliated structure with ambient variety $X$ and $\phi: X\dashrightarrow X'$ a birational map$/U$. Assume that $\phi$ is $K_{\Aa}$-negative, $\Aa':=\phi_*\Aa$, and let $X\rightarrow Z$ be a $K_{\Aa'}$-Mori fiber space$/U$. Then we say that $\Aa'\rightarrow Z$ is a \emph{Mori fiber space} of $\Aa/U$.
\end{defn}

We need the following two auxiliary results.

\begin{lem}[{\cite[Lemma 3.16(2)]{Cas+24}}]\label{lem: weak lc model only check codim 1}
 Let $\Aa/U$ be an adjoint foliated structure with ambient variety $X$ and let $\phi: X\dashrightarrow X'$ be a birational map$/U$ which does not extract any divisor. Let $\Aa':=\phi_*\Aa$. Assume that $K_{\Aa'}$ is nef$/U$. 

If $a(D,\Aa)<a(D,\Aa')$ for any prime divisor $D$ on $X$ that is exceptional$/X'$, then $\Aa'/U$ is a minimal model of $\Aa/U$.
\end{lem}

\begin{lem}[{cf. \cite[Lemma 5.4]{LMX25}}]\label{lem: gmmp scaling numbers go to 0}
Let $(X,B)/U$ be an lc pair. Let $H\geq 0$ be an $\Rr$-divisor on $X$ such that $(X,B+H)$ is lc and $K_X+B+H$ is nef$/U$. Assume that for any $\mu\in [0,1]$, either $(X,B+\mu H)/U$ has a log minimal model, or $K_X+B+\mu H$ is not pseudo-effective$/U$. Then there exists a $(K_X+B)$-MMP$/U$ with scaling of $H$ which terminates with a minimal model of $(X,B)/U$ after finitely many steps.
\end{lem}

\subsection{Locally stable families}

\begin{defn}[Slc]\label{defn: slc pair}
A \emph{semi-pair} $(X,B)$ consists of a demi-normal scheme $X$ and an $\Rr$-divisor $B\geq 0$ such that $K_X+B$ is $\Rr$-Cartier. Let $\nu: X^\nu\rightarrow X$ be the normalization of $X$ and let $D\subset X$, $D^\nu\subset X^\nu$ be the conductors. We say that $(X,B)$ is \emph{slc} if 
\begin{enumerate}
    \item $\Supp B$ does not contain any irreducible component of the conductor $D$.
    \item $(X^\nu,B^\nu+D^\nu)$ is lc, where $B^\nu$ is the divisorial part of $\nu^{-1}(B)$.
\end{enumerate}
We say that $(X,B)$ is \emph{stable} if $(X,B)$ is projective, slc, and $K_X+B$ is ample. 
\end{defn}

\begin{defn}[Mumford divisor]\label{defn: mumford divisor}
Let $d$ be a positive integer and let $f: X\rightarrow Z$ be a flat morphism over a reduced scheme such that the fibers of $f$ are reduced, connected, $S_2$, and of pure dimension $d$. A \emph{Mumford divisor}$/Z$ is a divisor $D$ on $X$ satisfying the following.
\begin{enumerate}
    \item (Equidimensionality) Every irreducible component of $\Supp D$ dominates an irreducible component of $Z$, and all nonempty fibers of the induced morphism $f|_D: \Supp D\rightarrow Z$ are of pure dimension $d-1$.
    \item (Mumford) $f$ is smooth near the generic point of $f^{-1}(z)\cap\Supp D$ for any point $z\in Z$.
    \item (Generic Cartier) $D$ is Cartier near the generic points of $f^{-1}(z)\cap\Supp D$ for any point $z\in Z$.
\end{enumerate}
An $\Rr$-divisor $B$ on $X$ is called a \emph{Mumford $\mathbb R$-divisor$/Z$} if $B=\sum a_iD_i$, where each $a_i\geq 0$ and each $D_i$ is a Mumford divisor$/Z$.
\end{defn}

\begin{defn}[Locally stable family]\label{defn: lsf}
A \emph{locally stable family} $f: (X,B)\rightarrow Z$ consists of a flat projective morphism $ f: X\rightarrow Z$ over a reduced scheme with demi-normal and connected fibers, and a Mumford $\Rr$-divisor$/Z$ $B$ on $X$, such that
\begin{enumerate}
\item $K_{X/Z}+B$ is $\Rr$-Cartier, and
\item $(X_z,B_z)$ is slc for any $z\in Z$, where $X_z=f^{-1}(z)$ and $B_z=B|_{X_z}$ (we refer the reader to \cite[Section 4.1]{Kol23} for the definition of such restriction).
\end{enumerate}
Here $K_{X/Z}$ is the relative canonical divisor corresponding to the relative dualizing sheaf $\omega_{X/Z}$ and we refer the reader to \cite[2.68]{Kol23} for more details.

A \emph{stable family} is a locally stable family $f: (X,B)\rightarrow Z$ such that $K_{X/Z}+B$ is ample$/Z$.
\end{defn}

\begin{defn}[Variation]
Let $f: (X,B)\rightarrow Z$ be a stable family. The \emph{variation} of $f: (X,B)\rightarrow Z$ is $\dim Z-d$ where $d$ is the dimension of a general isomorphism equivalence class of the log fibers. We say that $f: (X,B)\rightarrow Z$ is \emph{of maximal variation} if the variation of $f: (X,B)\rightarrow Z$ equals to $\dim Z$.
\end{defn}

\begin{defn}[Base change]\label{defn: base change}
    Let $f: (X,B)\rightarrow Z$ be a (locally) stable family. A \emph{base change} of $f: (X,B)\rightarrow Z$ is $f': (X',B')\rightarrow Z'$ satisfying the following:
    \begin{itemize}
        \item There is a morphism $h_Z: Z'\rightarrow Z$ such that $X'=X\times_ZZ'$, and $h: X'\rightarrow X$ is the induced morphism.
        \item $B'$ is the unique $\Rr$-divisor such that $K_{X'/Z'}+B'=h^*(K_{X/Z}+B)$.
    \end{itemize}
By definition (cf. \cite[Theorem 4.8]{Kol23}), $f': (X',B')\rightarrow Z'$ is also a (locally) stable family. If $h_Z$ is birational, then we say that $f': (X',B')\rightarrow Z'$ is a \emph{birational base change} of $f: (X,B)\rightarrow Z$.
\end{defn}

\begin{defn}[$(d,\Ii,v)$-stable family]
Let $d$ be a positive integer, $\Ii\subset [0,1]$ a set such that $\Ii=\Ii_+$, $v$ a positive real number, and $(X,B)$ a semi-pair. We say that $(X,B)$ is a \emph{$(d,\Ii,v)$-stable semi-pair} if  $\dim X=d$, $(X,B)$ is stable, the coefficients of $B$ belong to $\Ii$, and  $\vol(K_X+B)=v$. A \emph{family of $(d,\Ii,v)$-stable semi-pairs} is a locally stable family $f: (X,B)\rightarrow Z$, such that
 \begin{enumerate}
     \item $B=\sum b_iB_i$ where each $b_i\in\Ii$ and each $B_i$ is a Mumford divisor$/Z$, and
     \item $(X_z,B_z)$ is a $(d,\Ii,v)$-stable semi-pair for any point $z\in Z$, where $X_z=f^{-1}(z)$ and $B_z=B|_{X_z}$.
 \end{enumerate}
We denote by $\mathfrak{S}(d,\Ii,v)$ the set of all $(d,\Ii,v)$-stable semi-pairs. We denote by $\mathfrak{S}(d,\Ii,v)$ the moduli functor by setting
 $$\mathfrak{S}(d,\Ii,v)(Z)=\{\text{families of }(d,\Ii,v)\text{-stable semi-pairs over }Z\}$$
 where $Z$ is any reduced scheme.
\end{defn}

\begin{defthm}[{cf. \cite[Theorem 4.1]{Kol23}, \cite[Definition-Theorem 2.12]{HJLL24}}]\label{defthm: kol23 4.1}
Let $d$ be a positive integer, $\Ii\subset [0,1]$ a DCC set such that $\Ii=\Ii_+$, and $v$ a positive real number. Then the functor $\mathfrak{S}(d,\Ii,v)$ has a projective coarse moduli space $M(d,\Ii,v)$. For any family $f: (X,B)\rightarrow Z$ of $(d,\Ii,v)$-stable semi-pairs, we call the induced morphism $Z\rightarrow M(d,\Ii,v)$ the corresponding \emph{moduli map}.
\end{defthm}

\begin{thm}[{\cite[Corollary 6.19]{KP17}, \cite[Theorem 2.13]{HJLL24}}]\label{kp17 6.19}
Let $d$ be a positive integer, $\Ii\subset [0,1]$ a DCC set such that $\Ii=\Ii_+$, and $v$ a positive real number. Then there exists a reduced scheme $Z$ and $f\in\mathfrak{S}(d,\Ii,v)(Z)$, such that the corresponding moduli map $Z\rightarrow M(d,\Ii,v)$ is a finite cover.
\end{thm}

\subsection{Iitaka dimensions and abundance}

\begin{defn}\label{defn: iitaka dimension}
Let $X$ be a normal projective variety and $D$ an $\Rr$-divisor on $X$. The \emph{Iitaka dimension} $\kappa(D)$ (resp. \emph{numerical Iitaka dimension} $\nu(D)$) of $D$ is defined in the following way. If $|\lfloor mD\rfloor|\not=\emptyset$ for some positive integer $m$ (resp. $D$ is pseudo-effective), then
$$\kappa(D):=\max\left\{k\in\mathbb N\middle| \underset{m\rightarrow+\infty}{\lim\sup}\frac{\dim H^0(X,\lfloor mD\rfloor)}{m^k}>0\right\}$$
$$\left(\text{resp. }\nu(D):=\max\left\{k\in\mathbb N\middle| A\text{ is Cartier}, \underset{m\rightarrow+\infty}{\lim\sup}\frac{\dim H^0(X,\lfloor mD\rfloor+A)}{m^k}>0\right\}\right).$$
Otherwise, we let $\kappa(D):=-\infty$ (resp. $\nu(D):=-\infty$). If $|D|_{\mathbb R}\not=\emptyset$, then we define $\kappa_{\iota}(D):=\kappa(D')$ for some $D'\in |D|_{\mathbb R}$. Otherwise, we define $\kappa_{\iota}(D):=-\infty$. $\kappa_{\iota}$ is well-defined by \cite[Section 2]{Cho08}. We retain the notation $\kappa_{\iota}$ because it is invariant under $\mathbb R$-linear equivalence, whereas $\kappa$ for an $\mathbb R$-divisor is defined using the round-downs $\lfloor mD\rfloor$ and is a distinct invariant.

Let $\Ff$ be a foliation on $X$. We define $\kappa(X):=\kappa(K_X)$ (resp. $\kappa(\Ff):=\kappa(K_{\Ff})$) as the \emph{Kodaira dimension} of $X$ (resp. $\Ff$), and define $\nu(X):=\nu(K_X)$ (resp. $\nu(\Ff):=\nu(K_{\Ff})$) as the \emph{numerical dimension} of $X$ (resp. $\Ff$). We say that $X$ (resp. $\Ff$) is \emph{abundant} if $\kappa(X)=\nu(X)$ (resp. $\kappa(\Ff)=\nu(\Ff)$).

Let $\pi: X\rightarrow U$ be a projective morphism between normal quasi-projective varieties and $D$ an $\Rr$-divisor on $X$. Let $X\xrightarrow{f} T\rightarrow\pi(U)$ be the Stein factorization and let $F$ be a general fiber of $f$. We define $\kappa(X/U,D):=\kappa(D|_F),\kappa_{\iota}(X/U,D):=\kappa_{\iota}(D|_F)$, and $\nu(X/U,D):=\nu(D|_F)$ as the \emph{Iitaka dimension}, \emph{invariant Iitaka dimension}, and \emph{numerical Iitaka dimension} of $D$ over $U$ respectively.  We say that $D$ is \emph{abundant}$/U$ if $\kappa_{\iota}(X/U,D)=\nu(X/U,D)$. Here by convention, if $\dim F=0$, then we define $\kappa(D|_F):=\kappa_{\iota}(D|_F):=\nu(D|_F):=0$.

We refer the reader to \cite[Chapters II,V]{Nak04}, \cite[Section 2]{Cho08}, and \cite[Section 2]{HH20} for basic properties of Iitaka dimensions.
\end{defn}

We need the following result which characterizes nef and abundant divisors.
\begin{lem}[{\cite[Lemma 2.8]{Hu20}}]\label{lem: equivalent definition of abundance}
    Let $X$ be a normal projective variety and $D$ a nef$/U$ $\mathbb R$-divisor on $X$. Then the followings are equivalent.
\begin{enumerate}
    \item $D$ is abundant$/U$.
    \item There exists a birational morphism $g: Y\rightarrow X$ and a contraction$/U$ $h: Y\rightarrow Z$ such that $g^*D\sim_{\mathbb R}h^*H$ for some big$/U$ and nef$/U$ $\mathbb R$-divisor $H$ on $Z$.
    \item For a very general fiber $F$ of the Stein factorization of $X\rightarrow U$, either $\dim F=0$, or for every projective birational morphism $\rho:F'\rightarrow F$ with $F'$ normal, every prime divisor $P$ on $F'$, and every $\epsilon>0$, there exists
    $$0\leq D_{\rho,P,\epsilon}\sim_{\mathbb R}\rho^*(D|_F)
    \quad\text{such that}\quad \mult_PD_{\rho,P,\epsilon}<\epsilon.$$
\end{enumerate}
\end{lem}

Finally, we recall the following version of the abundance conjecture which is one assumption in many of our results.

\begin{conj}[Abundance]\label{conj: abundance}
    Let $(X,B)$ be a klt pair. Then $\kappa_{\iota}(K_X+B)=\nu(K_X+B)$.
\end{conj}

\subsection{Nakayama-Zariski decomposition}

\begin{defn}
    Let $\pi\colon X\rightarrow U$ be a projective morphism from a normal quasi-projective variety to a quasi-projective variety, $D$ a pseudo-effective$/U$ $\Rr$-Cartier $\Rr$-divisor on $X$, and $P$ a prime divisor on $X$. We define $\sigma_{P}(X/U,D)$ as in \cite[Definition 3.1]{LX25a} by considering $\sigma_{P}(X/U,D)$ as a number in  $[0,+\infty)\cup\{+\infty\}$. We define $N_{\sigma}(X/U,D)=\sum_Q\sigma_Q(X/U,D)Q$
    where the sum runs through all prime divisors on $X$ and consider it as a formal sum of divisors with coefficients in $[0,+\infty)\cup\{+\infty\}$. We say that $D$ is \emph{movable$/U$} if $N_{\sigma}(X/U,D)=0$. 
    
    If $U$ is a point then we denote by $\sigma_P(D):=\sigma(X/U,D)$ and $N_{\sigma}(D):=N_{\sigma}(X/U,D)$. By \cite[III 1.5 Lemma]{Nak04}, if $U$ is a point, then $N_{\sigma}(D)$ is an $\mathbb R$-divisor.
\end{defn}

We need the following properties on the negative part of Nakayama-Zariski decompositions.

\begin{lem}[{\cite[Lemma 2.25]{LMX25}}]\label{lem: nz for lc divisor}
Let $X\rightarrow U$ be a projective morphism from a normal quasi-projective variety to a quasi-projective variety and $\phi\colon X\dashrightarrow X'$ a birational map$/U$. Let $D$ be an $\Rr$-Cartier $\Rr$-divisor on $X$ such that $\phi$ is $D$-negative and $D':=\phi_\ast D$. Then:
\begin{enumerate}
    \item The divisors contracted by $\phi$ are contained in $\Supp N_{\sigma}(X/U,D)$.
    \item If $D'$ is movable$/U$, then $\Supp N_{\sigma}(X/U,D)$ is the set of all $\phi$-exceptional divisors.
    \end{enumerate}
\end{lem}

\begin{lem}[{cf. \cite[Lemma 3.21]{Cas+25a}}]\label{lem: if contract n then movable}
Let $X\rightarrow U$ be a projective morphism from a normal quasi-projective variety to a quasi-projective variety and $\phi: X\dashrightarrow X'$ a birational map$/U$ which does not extract any divisor. Let $D$ be an $\Rr$-Cartier $\Rr$-divisor on $X$ such that $D':=\phi_*D$ is $\Rr$-Cartier and $\phi$ contracts all divisors that are contained in  $\Supp N_{\sigma}(X/U,D)$. Then $D'$ is movable$/U$.
\end{lem}

\begin{lem}[{cf. \cite[Lemma 3.22]{Cas+25a}}]\label{lem: limit of nakayama-zariski decomposition}
    Let $X\rightarrow U$ be a projective morphism between normal quasi-projective varieties and let $C,D$ be two pseudo-effective$/U$ $\Rr$-Cartier $\Rr$-divisors on $X$. Assume that $N_{\sigma}(X/U,D)$ is an $\Rr$-divisor, i.e. it does not have $+\infty$ as a coefficient. Then there exists a positive real number $s_0$ and a reduced divisor $E$ satisfying the following.
    \begin{enumerate}
        \item $\Supp N_{\sigma}(X/U,C+sD)=E$ for any $0<s\leq s_0$.
        \item $\Supp N_{\sigma}(X/U,C)\subset E$.
    \end{enumerate}
\end{lem}

\section{Abundance of foliations}\label{sec: abu f}

The goal of this section is to prove Theorem \ref{thm: main nonalg I}.

\subsection{Abundant and semi-ampleness}

The following theorem is the real coefficient version of the recently proven Prokhorov-Shokurov-Bakker-Filipazzi-Mauri-Tsimerman's $\bb$-semi-ampleness theorem \cite[Theorem 1.5]{BFMT25}.
\begin{thm}\label{thm: real coefficient bsa}
    Let $(X,B)$ be a projective lc pair and let $f: X\rightarrow Z$ be a contraction such that $K_X+B\sim_{\mathbb R,Z}0$. Let $(Z,B_Z,\Mm)$ be the generalized pair induced by the canonical bundle formula of $f: (X,B)\rightarrow Z$,
    $$K_X+B\sim_{\mathbb R}f^*(K_Z+B_Z+\Mm_Z).$$
    Then $\Mm$ is $\bb$-semi-ample.
\end{thm}
\begin{proof}
Write $B=\sum_{i=1}^m v_i^0B_i$ where $B_i\geq 0$ are Weil divisors. Let $\bm{v}_0:=(v_1^0,\dots,v_m^0)$, and let $V\ni\bm{v}_0$ be the rational envelope of $\bm{v}_0$ in $\mathbb R^m$. Let $B(\bm{v}):=\sum_{i=1}^mv_iB_i$ for any $\bm{v}\in\mathbb R^m$. By \cite[Theorem 3.3]{HLX25} and \cite[Theorem 5.6]{HLS26}, there exists an open subset $U\ni\bm{v}_0$ of $V$, such that for any $\bm{v}\in U$,
\begin{itemize}
\item $f: (X,B(\bm{v}))\rightarrow Z$ is an lc-trivial fibration and $(X,B(\bm{v}))$ is lc, and
\item Let $(Z,B_Z(\bm{v}),\Mm(\bm{v}))$ be the generalized pair induced by the canonical bundle formula of $f: (X,B(\bm{v}))\rightarrow Z$,
$$K_X+B(\bm{v})\sim_{\mathbb R}f^*(K_Z+B_Z(\bm{v})+\Mm(\bm{v})_Z),$$
then for any vectors $\bm{v}_1,\dots,\bm{v}_k\in U$ and positive real numbers $a_1,\dots,a_k$ such that $\sum_{i=1}^ka_i=1$, we have
$$\Mm\left(\sum_{i=1}^ka_i\bm{v}_i\right)=\sum_{i=1}^ka_i\Mm(\bm{v}_i).$$
\end{itemize}
Pick $\bm{v}_i^0\in\mathbb Q^m\cap U$ and $a_i^0$ such that $\sum_{i=1}^ka_i^0\bm{v}_i^0=\bm{v}_0$ and $\sum_{i=1}^ka_i^0=1$. By \cite[Theorem 1.5]{BFMT25}, $\Mm(\bm{v}_i^0)$ is $\bb$-semi-ample for any $i$, hence
$$\Mm=\Mm\left(\sum_{i=1}^ka_i^0\bm{v}_i^0\right)=\sum_{i=1}^ka_i^0\Mm(\bm{v}_i^0)$$
is $\bb$-semi-ample.
\end{proof}

\subsection{Abundance for stable families}

The following proposition is the real coefficient version of \cite[Corollary 6.19]{KP17}.

\begin{prop}\label{prop: base change real coefficient stable family}
Let $f: (X,B)\rightarrow Z$ be a stable family such that $Z$ is normal. Then there exist two stable families $f': (X',B')\rightarrow Z'$, $f'': (X'',B'')\rightarrow Z''$, a finite surjective morphism $g: Z''\rightarrow Z$, and a dominant morphism $g': Z''\rightarrow Z'$, such that
\begin{enumerate}
    \item $f'': (X'',B'')\rightarrow Z''$ is the base change of both $f: (X,B)\rightarrow Z$ and $f': (X',B')\rightarrow Z'$,
\item $Z'$ and $Z''$ are normal,
    \item $f': (X',B')\rightarrow Z'$ is of maximal variation, and
    \item if the general fiber of $f$ is normal, then the general fibers of $f'$ and $f''$ are normal.
\end{enumerate}
\end{prop}
\begin{proof}
Let $F$ be a general fiber of $f$ and let $B_F:=B|_F$. If $\dim F=0$, then $f$ is finite \'etale and $B=0$; a finite \'etale cover trivializes this family of reduced points, and we take $Z'$ to be a point. Assume henceforth that $d:=\dim F>0$. Let $\Ii'$ be the coefficient set of $B_F$, let $\Ii:=\Ii'_+$, and let $v:=\vol(K_F+B_F)$. Let $Z\rightarrow M(d,\Ii,v)$ be the moduli map of $f$. By Theorem \ref{kp17 6.19}, there is a stable family $\ff:(\mathcal X,\mathcal B)\rightarrow\mathcal Z$ whose moduli map $\mathcal Z\rightarrow M(d,\Ii,v)$ is finite and surjective.

Let $T$ be the normalization of a component of $Z\times_{M(d,\Ii,v)}\mathcal Z$ dominating $Z$. The two pulled-back stable families on $T$ have the same moduli map. By \cite[Proposition 2.15]{HJLL24}, after a further finite surjective cover $Z''\rightarrow T$, followed by normalization, these families are isomorphic. Fix this isomorphism and denote the common family by $(X'',B'')/Z''$.

Let $Z'$ be the normalization of the closure of the image of $Z''\rightarrow\mathcal Z$, and let $(X',B')/Z'$ be the induced stable family. The morphism $g:Z''\rightarrow Z$ is finite and surjective, and the induced morphism $g':Z''\rightarrow Z'$ is dominant. Its image is closed, and hence $g'$ is surjective, whenever $Z$ is projective, as in the application below. Since $Z'\rightarrow M(d,\Ii,v)$ is finite onto its image, the family on $Z'$ has maximal variation. The chosen isomorphism gives the two required base-change identifications. Normality of the general fibers is preserved under these base changes in characteristic zero.
\end{proof}

\begin{lem}\label{lem: abudant stable family}
    Let $f: (X,B)\rightarrow Z$ be a stable family such that $Z$ is projective and the general fiber of $f$ is normal. Then $K_{X/Z}+B$ is nef and abundant.
\end{lem}
\begin{proof}
By Proposition \ref{prop: base change real coefficient stable family}, there exists a stable family of maximal variation $f': (X',B')\rightarrow Z'$, and two projective surjective morphisms $h: X''\rightarrow X$, $h': X''\rightarrow X'$, such that the general fiber of $f'$ is normal and 
$$h^*(K_{X/Z}+B)=h'^*(K_{X'/Z'}+B').$$
Here $X'$ is the maximal-variation family to which the positivity result is applied, whereas $X''$ is the common base change on which the two pullbacks are compared. Since the invariant Iitaka dimension and the numerical dimension are unchanged under projective surjective pullback, the pullback identity gives the equalities below.
If the fibers have dimension zero, then $f$ is finite \'etale and $B=0$, so $K_{X/Z}+B=0$ and the assertion is immediate. Otherwise, the normal generic log fiber $(X'_{\eta},B_{\eta}')$ is lc. Since $Z$ is projective, $Z'$ is projective. By \cite[Proposition 2.18]{HJLL24} (rational coefficient case \cite[Proposition 2.15]{PX17}), $K_{X'/Z'}+B'$ is big and nef. Its pullback to $X''$ is nef, and nefness descends under the surjective morphism $h$. Thus $K_{X/Z}+B$ is nef. Moreover,
$$\kappa_{\iota}(K_{X/Z}+B)=\kappa_{\iota}(K_{X'/Z'}+B')=\nu(K_{X'/Z'}+B')=\nu(K_{X/Z}+B).$$
\end{proof}

\subsection{Transformation of generalized foliated quadruples}

The following proposition is the generalized foliated quadruple version of \cite[Proposition 3.16]{HJLL24} which characterizes generalized foliated quadruples under finite covers.

\begin{prop}\label{prop: cover formula}
    Let $\Aa:=(X,\Ff,B,\Mm)$ be a generalized foliated quadruple and let $h: X'\rightarrow X$ be a finite cover between normal varieties. Let $\Ff':=h^{-1}\Ff$ and $\Mm':=h^*\Mm$. For any prime divisor $D$ on $X$, let $r_D$ be the ramification index of $h$ along $D$. Suppose that any codimension one component of the branch locus of $h$ is $\Ff$-invariant. Let $K_{\Ff'}+B'+\Mm'_{X'}:=h^*(K_{\Ff}+B+\Mm_X)$ and $\Aa':=(X',\Ff',B',\Mm')$. Then:
    \begin{enumerate}
         \item If $\Aa$ is lc (resp. klt, canonical), then $\Aa'$ is also lc (resp. klt, canonical).
        \item If $\Aa'$ is lc (resp. klt), then $\Aa$ is lc (resp. klt).
    \end{enumerate}
\end{prop}
\begin{proof}
Let $g: Y\rightarrow X$ be a birational morphism and $E$ a prime divisor on $Y$. Let $Y'$ be the normalization of the main component of $Y\times_{X}X'$ associated with $g': Y'\rightarrow X'$ and $h': Y'\rightarrow Y$. Let $E'$ be any prime divisor on $Y'$ that dominates $E$. Such prime divisor exists as $h'$ is finite. Let $\Ff_Y:=g^{-1}\Ff$ and $\Ff_{Y'}:=g'^{-1}\Ff'$.

Let $r_{E'}$ be the ramification index of $h'$ along $E'$. By \cite[Proposition 3.15(1)]{HJLL24}, near the generic point of $E'$, we have
\begin{align*}
K_{\Ff_{Y'}}&=g'^*K_{\Aa'}+a(E',\Aa')E'=g'^*h^*K_{\Aa}+a(E',\Aa')E'\\
  &=h'^*g^*K_{\Aa}+a(E',\Aa')E'=h'^*(K_{\Ff_Y}-a(E,\Aa)E)+a(E',\Aa')E'\\
  &=K_{\Ff_{Y'}}-\epsilon_{\Ff'}(E')(r_{E'}-1)E'-a(E,\Aa)r_{E'}E'+a(E',\Aa')E'.
\end{align*}
So
$$a(E',\Aa')+\epsilon_{\Ff'}(E')=r_{E'}(a(E,\Aa)+\epsilon_{\Ff'}(E'))=r_{E'}(a(E,\Aa)+\epsilon_{\Ff}(E)).$$
Therefore, $a(E,\Aa)\geq-\epsilon_{\Ff}(E)$ (resp. $>-\epsilon_{\Ff}(E)$) if and only if $a(E',\Aa')\geq-\epsilon_{\Ff'}(E')$ (resp. $>-\epsilon_{\Ff'}(E')$), and if $a(E,\Aa)\geq 0$ (resp. $>0$), then $a(E',\Aa')\geq 0$ (resp. $>0$). Moreover, $E$ is exceptional$/X$ if and only if $E'$ is exceptional$/X'$. By \cite[Lemma 3.6]{HJLL24}, we get (2). By \cite[Lemma 2.22]{Kol13}, for any prime divisor $E'$ over $X'$, we may find $Y',Y,g,g',h'$ and $E$ as above. (1) follows.
\end{proof}

The following proposition is the generalized foliated quadruple version of \cite[Proposition 3.16]{HJLL24} which concerns the behavior of generalized foliated quadruples under birational base changes.

\begin{lem}\label{lem: ACSS exist core mod equidim}
    Let $\Aa/U:=(X,\Ff,B,\Mm)/U$ be an lc generalized foliated quadruple. Then 
    \begin{enumerate}
        \item $\Aa$ has a $\Qq$-factorial ACSS modification $h: (\Aa';G)/Z \to \Aa$.
        \item Let $X'$ be the ambient variety of $\Aa'$. Then there exists a core modification $\bar h: (\bar\Aa;\bar G)/Z\to \Aa$ associated with $\bar f: \bar X\to Z$ such that $\bar X$ is the core model of $(h,f: X'\to Z)$ associates with $(\bar h, \bar f)$ and $\bar f$ is equidimensional.
    \end{enumerate}
\end{lem}
\begin{proof}
(1) follows from \cite[Theorem 8.2.2]{CHLX23}. (2) By \cite[Lemma A.22]{LMX25}, the core model of $(h,f:X'\to Z)$ associated with $(\bar h,\bar f)$ induces a core modification $(\bar\Aa;\bar G)/Z\rightarrow\Aa$. Since $f: X'\to Z$ is equidimensional, the induced contraction $\bar f: \bar X \to Z$ is also equidimensional. (2) follows.
\end{proof}

\subsection{Proof of Theorem \ref{thm: main nonalg I}}

The following proposition is the key of the proof of Theorem \ref{thm: main nonalg I}. Some ideas of the proof of Proposition \ref{prop: nef relative sa imply abundant} originated from \cite[Proof of Theorem 4.6]{FS25} although the details have many differences.

\begin{prop}\label{prop: nef relative sa imply abundant}
Let $(X,\Ff,B)$ be a projective lc algebraically
integrable foliated triple. Assume that $(X,\Ff,B;G)/Z$ satisfies Property $(*)$ associated with
an equidimensional contraction $f: X\rightarrow Z$ and $K_{\Ff}+B$ is semi-ample$/Z$. Then $K_{\Ff}+B$ is abundant.
\end{prop}
\begin{proof}
\noindent\textbf{Step 1.} In this step we construct a stable family after base change.

By \cite[Proposition 3.16]{HJLL24} (rational coefficient case \cite[Proposition 3.5]{FS25}), there is a finite cover $\tau:X_1\rightarrow X$ such that, after writing
$$K_{\Ff_1}+B_1=\tau^*(K_{\Ff}+B),$$
the induced morphism $f_1:(X_1,B_1)\rightarrow Z_1$ is a locally stable family. Both $\kappa_{\iota}$ and $\nu$ are invariant under finite surjective pullback, so $K_{\Ff_1}+B_1$ is abundant if and only if $K_{\Ff}+B$ is abundant. After this finite cover, take a projective resolution of the base and form the induced locally stable base change. Relative log canonical divisors commute with this base change (Definition \ref{defn: base change}); the new total space is normal because the base and the generic fiber are normal (cf.\ \cite[10.78]{Kol23}). Invariant Iitaka dimension and numerical dimension are unchanged under this projective surjective pullback. After relabelling, we may therefore assume that $f:(X,B)\rightarrow Z$ is a locally stable family and $Z$ is smooth and projective. We let $g: X\rightarrow Y$ be the ample model$/Z$ of $K_{X/Z}+B$, then we have
$$K_X+B\sim_{\mathbb R,Y}K_{X/Z}+B\sim_{\mathbb R,Y}0.$$
Let $h:Y\rightarrow Z$ be the induced contraction and let $\Ff_Y$ be the foliation induced by $h$. We let $(Y,B_Y,\Mm)$ be the lc generalized pair induced by the canonical bundle formula of $g:(X,B)\rightarrow Y$,
$$K_X+B\sim_{\mathbb R}g^*(K_Y+B_Y+\Mm_Y).$$

For any closed point $z\in Z$, let $H_1,\dots,H_{\dim Z}$ be general hyperplane sections on $Z$ such that $z\in H_i$ for each $i$ and let $H_z:=\sum_{i=1}^{\dim Z} H_i$. Since $f: (X,B)\rightarrow Z$ is a locally stable family, $X$ is normal, and $Z$ is smooth, $(X,B+f^*H_z)$ is lc. Since $(Y,B_Y+h^*H_z,\Mm)$ is the generalized pair induced by the canonical bundle formula
$$K_X+B+f^*H_z\sim_{\mathbb R}g^*(K_Y+B_Y+h^*H_z+\Mm_Y),$$
$(Y,B_Y+h^*H_z,\Mm)$ is lc. By Theorem \ref{thm: real coefficient bsa}, $\Mm$ is $\bb$-semi-ample. Let $\alpha:T\rightarrow Y$ be a birational morphism such that $\Mm$ descends to $T$. Then there exists $0\leq L_{z,T}\sim_{\mathbb R}\Mm_T$ such that $(Y,B_Y+h^*H_z+L_z)$ is lc, where $L_z:=\alpha_*L_{z,T}$. Let $Y_z:=h^{-1}(z)$, $B_{Y_z}:=B_Y|_{Y_z}$, and $L_{Y_z}:=L_z|_{Y_z}$. By adjunction, $(Y_z,B_{Y_z}+L_{Y_z})$ is slc. By \cite[Corollary 4.45]{Kol23}, there exists an open subset $U_z\ni z$ with $V_z:=h^{-1}(U_z)$, such that $h|_{V_z}: (Y,B_Y+L_z)|_{V_z}\rightarrow U_z$ is a locally stable family, hence a stable family.

By Noetherian property, there exist finitely many closed points $z_1,\dots,z_k\in Z$ such that $\{U_{z_i}\}$ is a cover of $Z$. Put
\[
D:=K_{Y/Z}+B_Y+\Mm_Y.
\]
If $\dim Y=\dim Z$, each $V_{z_i}\to U_{z_i}$ is finite flat of degree one, hence $h$ is an isomorphism. A family of relative dimension zero has zero Mumford boundary, so $N_i:=B_Y+L_{z_i}\geq0$ satisfies $N_i\sim_{\mathbb R}D$ and $N_i|_{V_{z_i}}=0$. Choosing $V_{z_i}$ to contain the generic point of any curve, or of the center of any prime divisor over $Y$, shows respectively that $D$ is nef and satisfies Lemma~\ref{lem: equivalent definition of abundance}(3). Thus $D$, and its pullback under $g$, are abundant.

Assume henceforth that $\dim Y>\dim Z$. Fix an irreducible closed subset $S\subset Y$ and choose $i$ such that $V:=V_{z_i}$ contains its generic point. Let $U:=U_{z_i}$, $L:=L_{z_i}$, $L_V:=L|_V$, and $B_V:=B_Y|_V$. Then $h_V:=h|_V:(V,B_V+L_V)\rightarrow U$ is a stable family with normal general fiber.

Let $F$ be a general fiber of $h$, let $d:=\dim F$, let $\Ii'$ be the set of coefficients of $B_V+L_V$, let $\Ii:=\Ii'_+$, and let $v:=\vol((K_Y+B_Y+\Mm_Y)|_F)$. Then $h_V: (V,B_V+L_V)\rightarrow U$ is a family of $(d,\Ii,v)$-stable semi-pairs. Thus there exists a moduli map $\phi_0: U\rightarrow M(d,\Ii,v)$ which naturally extends to a rational map $\phi: Z\dashrightarrow M(d,\Ii,v)$. By Theorem \ref{kp17 6.19}, there exists a reduced scheme $\mathcal{Z}$ and a family $\ff: (\mathcal{X},\mathcal{D})\rightarrow\mathcal{Z}$ such that $\mathfrak{f}\in\mathfrak{S}(d,\Ii,v)(\mathcal{Z})$ and the induced moduli map $\mathcal{Z}\rightarrow M(d,\Ii,v)$ is a finite cover. Let $T_U$ be the normalization of a component of $U\times_{M(d,\Ii,v)}\mathcal Z$ dominating $U$.
The two stable families on $T_U$ have the same moduli map. By \cite[Proposition 2.15]{HJLL24}, there is a further finite surjective cover $U'\rightarrow T_U$, with $U'$ normal, on which these families are isomorphic. Fix this isomorphism.
Compactify $U'$ and take the normalization of the closure of the graph of its maps to $Z$ and $\mathcal Z$. This gives a normal projective variety $Z'$ containing $U'$ and extending both maps. Let $h':(Y',D')\rightarrow Z'$ be the stable family obtained from $(\mathcal X,\mathcal D)/\mathcal Z$ by base change. Its restriction to $U'$ is, by the chosen isomorphism, the base change of $(V,B_V+L_V)/U$.

Let $\pi_Z: Z'\rightarrow Z$ be the associated morphism. Then $\pi_Z$ is generically finite, and we may let $Z'\rightarrow Z''\rightarrow Z$ be the Stein factorization of $\pi_Z$. Let $Y'':=Y\times_{Z}Z''$, $\widetilde Y=Y\times_{Z}Z'$, and $\tau: Y''\rightarrow Y$, $\xi: \widetilde Y\rightarrow Y''$, $f'': Y''\rightarrow Z''$, and $\widetilde f: \widetilde Y\rightarrow Z'$ be the associated morphisms. Write $\Mm':=\tau^*\Mm$.
Since $U'\rightarrow U$ is finite, its open immersion into $\pi_Z^{-1}(U)$ is also proper, hence closed. The space $\pi_Z^{-1}(U)$ is irreducible, so $\pi_Z^{-1}(U)=U'$.
Let $\pi:=\tau\circ\xi:\widetilde Y\rightarrow Y$, and let $\Ff''$ and $\widetilde\Ff$ be the foliations induced by $f''$ and $\widetilde f$. The base-change maps fit into the diagram
\[
\begin{tikzcd}
\widetilde Y \arrow[r,"\xi"] \arrow[d,"\widetilde f"'] &
Y'' \arrow[r,"\tau"] \arrow[d,"f''"'] & Y \arrow[d,"h"] \\
Z' \arrow[r] & Z'' \arrow[r] & Z .
\end{tikzcd}
\]
Write
$$K_{\Ff''}+B''+\Mm'_{Y''}=\tau^*(K_{\Ff_Y}+B_Y+\Mm_Y),K_{\widetilde{\Ff}}+\widetilde B+\Mm'_{\widetilde{Y}}=\xi^*(K_{\Ff''}+B''+\Mm'_{Y''}).$$

Let $\widetilde T$ be the normalization of the main component of $T\times_Y\widetilde Y$, with projections $\pi_T:\widetilde T\to T$ and $\widetilde\alpha:\widetilde T\to\widetilde Y$.
We claim that $\widetilde B\geq0$. Fix a prime divisor $Q$ on $\widetilde Y$ and let $Q_T$ be its strict transform on $\widetilde T$. Choose a closed point $z$ of the closure of $(h\circ\pi)(Q)$. Since $\Mm_T$ is semi-ample, in the construction above we may also choose $L_{z,T}$ to avoid the generic point of $\pi_T(Q_T)$. Set $\widetilde L_z:=\widetilde\alpha_*\pi_T^*L_{z,T}$, so $\mult_Q\widetilde L_z=0$.

The open set $U_z$ contains the generic point of the closure of $(h\circ\pi)(Q)$. Locally stable base change gives an effective boundary $\widetilde D_z$ near the generic point of $Q$. Since the fibers are reduced and $\Mm_T-L_{z,T}$ is $\mathbb R$-principal, we have there
\[
\begin{aligned}
K_{\widetilde\Ff}+\widetilde D_z&=\pi^*(K_{\Ff_Y}+B_Y+L_z),\\
\pi^*(\Mm_Y-L_z)&=\Mm'_{\widetilde Y}-\widetilde L_z.
\end{aligned}
\]
Comparing with the defining pullback identity gives $\widetilde D_z=\widetilde B+\widetilde L_z$, hence $\mult_Q\widetilde B=\mult_Q\widetilde D_z\geq0$. This proves the claim.

Let $\widetilde L:=\widetilde\alpha_*\pi_T^*L_{z_i,T}\geq0$ and $\widetilde D:=\widetilde B+\widetilde L$. The chosen family isomorphism induces a birational map $\psi:\widetilde Y\dashrightarrow Y'$ that identifies $(\widetilde Y,\widetilde D)$ and $(Y',D')$ over $U'$.

\medskip

\noindent\textbf{Step 2.} In this step we compare the two stable families.

We let $p: W\rightarrow \widetilde Y$ and $q: W\rightarrow Y'$ be a resolution of indeterminacy of $\psi$. Thus
\[
\begin{tikzcd}
& W \arrow[dl,"p"'] \arrow[dr,"q"] & \\
\widetilde Y \arrow[rr,dashed,"\psi"'] \arrow[dr,"\widetilde f"'] &&
Y' \arrow[dl,"h'"] \\
& Z' &
\end{tikzcd}
\]
is a common resolution. Let $\Ff'$ be the foliation induced by $h'$. Then $\psi_*\widetilde{\Ff}=\Ff'$. Write
$$p^*(K_{\widetilde{\Ff}}+\widetilde D)=q^*(K_{\Ff'}+D')+E.$$
For any prime divisor $Q$ on $\widetilde Y$, if $Q$ is not $\widetilde\Ff$-invariant, then 
$$a(Q,\widetilde\Ff,\widetilde D)=a(Q,\Ff',D')$$
as $\psi$ is an isomorphism over the generic point of $Q$. If $Q$ is $\widetilde\Ff$-invariant, then
$$a(Q,\widetilde\Ff,\widetilde{D})=-\mult_Q\widetilde D\leq 0\leq a(Q,\Ff',D')$$
as $(Y',\Ff',D')$ is lc, which in turn is because $h': (Y',D')\rightarrow Z'$ is stable family. Therefore,
$$p_*E=\sum_{Q\subset\widetilde Y}\left(a(Q,\Ff',D')-a(Q,\widetilde\Ff,\widetilde D)\right)Q\geq 0.$$
The morphisms $p$ and $q$ are over $Z'$, and $K_{Y'/Z'}+D'$ is ample$/Z'$. Hence $q^*(K_{Y'/Z'}+D')$ is nef over $\widetilde Y$, and
$$E\sim_{\mathbb R,\widetilde Y}-q^*(K_{\Ff'}+D')=-q^*(K_{Y'/Z'}+D')$$
is anti-nef over $\widetilde Y$. By the negativity lemma, $E\geq 0$. Moreover, $p_*E$ is vertical$/Z'$.

Write $r:=\pi\circ p$ and $H:=q^*(K_{Y'/Z'}+D')$. By Lemma~\ref{lem: abudant stable family}, $H$ is nef and abundant, and
\[
r^*D\sim_{\mathbb R}H+E,\qquad E\geq0,\qquad E|_{r^{-1}(V)}=0.
\]
The last equality follows from the isomorphism of pairs over $U'$ and $\pi^{-1}(V)=\widetilde f^{-1}(U')$, with compatible canonical divisors; it includes all $p$-exceptional components of $E$.

\medskip

\noindent\textbf{Step 3.} In this step we prove that $D$ is nef and abundant.

For an integral curve $C\subset Y$, carry out the construction with $S=C$. Choose a curve $\Gamma\subset W$ dominating $C$ and put $e:=\deg(\Gamma/C)$. Since $\Gamma$ meets $r^{-1}(V)$, it is not contained in $\Supp E$, and
\[
e(D\cdot C)=H\cdot\Gamma+E\cdot\Gamma\geq0.
\]
Since $C$ is arbitrary, $D$ is nef. To prove abundance, let $\mu_0:Y_0\to Y$ be any projective birational morphism with $Y_0$ normal, and let $P$ be a prime divisor on $Y_0$. Apply the construction in \textbf{Steps 1--2} with $S:=\mu_0(P)$. A resolution $W_0$ of the main component of $W\times_Y Y_0$ has projective morphisms $a:W_0\to W$ and $b:W_0\to Y_0$, where $a$ is birational, $b$ is generically finite, and $r\circ a=\mu_0\circ b$. Put $H_0:=a^*H$, $E_0:=a^*E$, and $m:=\deg b$. Then $b^*\mu_0^*D\sim_{\mathbb R}H_0+E_0$. Since $E|_{r^{-1}(V)}=0$ and $V$ contains the generic point of $S$, the coefficient of $E_0$ along every prime divisor dominating $P$ is zero.

By Lemma~\ref{lem: equivalent definition of abundance}(2), after replacing $W_0$ by a higher model we have $H_0\sim_{\mathbb R}t^*J$ for a contraction $t:W_0\to T_0$ to a normal projective variety and a big and nef $\mathbb R$-divisor $J$ on $T_0$. Choose $J\sim_{\mathbb R}A+C_0$ with $A$ ample and $C_0\geq0$. Since
\[
J-\frac1nC_0\sim_{\mathbb R}\left(1-\frac1n\right)J+\frac1nA
\]
is ample for every $n\geq1$, we may write $H_0\sim_{\mathbb R}A_n+\frac1nC$, where $A_n$ is semi-ample and $C:=t^*C_0\geq0$. Choose $0\leq H_n\sim_{\mathbb R}A_n$ containing none of the finitely many prime divisors dominating $P$. The projection and norm formulas give
\[
\begin{gathered}
N_{0,n}:=\frac1m b_*\left(H_n+\frac1nC+E_0\right)\geq0,
\qquad N_{0,n}\sim_{\mathbb R}\mu_0^*D,\\
\mult_PN_{0,n}=\frac1{mn}\mult_P(b_*C)\longrightarrow0.
\end{gathered}
\]
Since $\mu_0$ and $P$ are arbitrary, Lemma~\ref{lem: equivalent definition of abundance}(3) implies that $D$ is abundant. Thus $K_{\Ff}+B\sim_{\mathbb R}g^*D$ is abundant, as required.
\end{proof}

\begin{thm}\label{thm: log abundance aif}
Assume Conjecture \ref{conj: abundance} in dimension $d$. Then for any projective lc algebraically integrable foliated triple $(X,\Ff,B)$ such that $\rk\Ff\leq d$, $K_{\Ff}+B$ is abundant. Moreover, if $X$ is potentially klt, then:
\begin{enumerate}
    \item $(X,\Ff,B)$ has a $\mathbb Q$-factorial minimal model or a Mori fiber space.
    \item We may run a $(K_{\Ff}+B)$-MMP with scaling of an ample $\mathbb R$-divisor and any such MMP terminates with either a minimal model or a Mori fiber space of $(X,\Ff,B)$. 
\end{enumerate}
\end{thm}
\begin{proof}
If $K_{\Ff}+B$ is not pseudo-effective, then
$$
\kappa_{\iota}(K_{\Ff}+B)=\nu(K_{\Ff}+B)=-\infty,
$$
so $K_{\Ff}+B$ is abundant. Moreover, when $X$ is potentially klt, \cite[Theorem 2.1.5]{Cas+25a} shows that any $(K_{\Ff}+B)$-MMP with scaling of an ample $\mathbb R$-divisor terminates with a Mori fiber space. Thus both assertions in the moreover part hold in this case.

We may therefore assume that $K_{\Ff}+B$ is pseudo-effective. Let $h: (X',\Ff',B';G)/Z\rightarrow (X,\Ff,B)$ be a $\mathbb Q$-factorial ACSS modification of $(X,\Ff,B)$ associated with a contraction $f: X'\rightarrow Z$. By \cite[Theorem 1.4]{HH20}, $(X',B'+G)/Z$ has a good minimal model. We run a $(K_{\Ff'}+B')$-MMP with scaling of an ample divisor. By \cite[Lemma 9.1.3(4)]{CHLX23}, this is also a sequence of steps of a $(K_{X'}+B'+G)$-MMP$/Z$. By \cite[Theorem 1.9]{Bir12}, this MMP terminates with a good minimal model $(X'',B''+G'')$ of $(X',B'+G)/Z$, where $B''$ and $G''$ denote the corresponding birational transforms. Let $\Ff''$ be the foliation induced by $X''\rightarrow Z$. By \cite[Lemma 9.1.3(3)]{CHLX23}, $(X'',\Ff'',B'';G'')/Z$ is $\mathbb Q$-factorial ACSS. By \cite[Proposition 3.6]{ACSS21}, $K_{\Ff''}+B''$ is semi-ample$/Z$. Proposition \ref{prop: nef relative sa imply abundant} implies that $K_{\Ff''}+B''$ is abundant, and hence so is $K_{\Ff}+B$. Moreover, $K_{\Ff''}+B''$ is nef: otherwise, \cite[Lemma 9.1.3(1)]{CHLX23} would give a negative extremal ray over $Z$, contradicting the relative semiampleness.

It remains to prove the moreover part in the pseudo-effective case. The model $(X'',\Ff'',B'')$ constructed above is a bs-minimal model of $(X,\Ff,B)$. Therefore, \cite[Theorem 2.3.1(4)]{Cas+25a} proves (2). For (1), let $\pi:\widetilde X\rightarrow X$ be a small $\mathbb Q$-factorialization, set $\widetilde\Ff:=\pi^{-1}\Ff$, and write $\widetilde B:=\pi^{-1}_*B$. Since $\pi$ is small, the model above is also a bs-minimal model of $(\widetilde X,\widetilde\Ff,\widetilde B)$. Applying \cite[Theorem 2.3.1(4)]{Cas+25a} to this triple gives a terminating MMP whose output is a $\mathbb Q$-factorial minimal model of $(X,\Ff,B)$. This proves (1).
\end{proof}

\begin{proof}[Proof of Theorem \ref{thm: main nonalg I}]
It is a special case of Theorem \ref{thm: log abundance aif}.
\end{proof}

\begin{proof}[Proof of Corollary \ref{cor: rank 3 abudnance}]
By \cite[1.1 Theorem]{KMM94}, Conjecture \ref{conj: abundance} holds in dimension $\leq 3$. (Note that \cite[1.1 Theorem]{KMM94} is stated for pairs with rational coefficients but the general case follows easily, cf. \cite[Theorem 2.20]{LX25b}). The corollary follows from Theorem \ref{thm: log abundance aif}.
\end{proof}

\section{Abundance for algebraically integrable adjoint foliated structures}\label{sec: abu aiafs}

The goal of this section is to prove Theorem \ref{thm: main abundance aiafs}.

\subsection{Abundant divisors} We first state two useful lemmas on klt varieties with polarization of abundant $\mathbb R$-divisors.

\begin{lem}\label{lem: abundant imply gmm}
    Let $(X,B)/U$ be a klt pair such that $K_X+B$ is abundant$/U$.
    \begin{enumerate}
        \item If $K_X+B$ is not pseudo-effective$/U$, then we may run a $(K_X+B)$-MMP$/U$ with scaling of an ample$/U$ $\mathbb R$-divisor and any such MMP terminates with a Mori fiber space of $(X,B)/U$.
        \item If $K_X+B$ is pseudo-effective$/U$, then we may run a $(K_X+B)$-MMP$/U$ with scaling of an ample$/U$ $\mathbb R$-divisor and any such MMP terminates with a good minimal model of $(X,B)/U$. Moreover, $(X,B)/U$ has a $\mathbb Q$-factorial good minimal model.
    \end{enumerate}
\end{lem}
\begin{proof}
Let $h: X'\rightarrow X$ be a small $\mathbb Q$-factorialization of $X$ and write $K_{X'}+B':=h^*(K_X+B)$. Let $A$ be an ample$/U$ $\mathbb R$-divisor on $X$ and let $A':=h^*A$. By \cite[Remark 2.9]{Bir12}, any sequence of steps of a $(K_X+B)$-MMP$/U$ with scaling of $A$ induces a $(K_{X'}+B')$-MMP$/U$ with scaling of $A'$. Let $\lambda_i$ be the scaling numbers of the $(K_{X'}+B')$-MMP$/U$ with scaling of $A'$. 

Since $A'$ is big, if $K_X+B$ is not pseudo-effective$/U$, then $K_{X'}+B'$ is not pseudo-effective$/U$. By \cite[Corollary 1.4.2]{BCHM10}, this MMP terminates, so that $(K_X+B)$-MMP$/U$ with scaling of $A$ terminates with a Mori fiber space of $(X,B)/U$. Thus we may assume that $K_X+B$ is pseudo-effective$/U$. By \cite[Lemma 2.13]{HH20} (rational coefficient case \cite[Theorem 4.3]{GL13}), $(X,B)$ has a good minimal model $(X_{\min},B_{\min})/U$. Since $h$ is small,  $(X_{\min},B_{\min})/U$ is also a good minimal model of $(X',B')/U$. By \cite[Theorem 1.9]{Bir12} and \cite[Corollary 1.4.2]{BCHM10}, any $(K_{X'}+B')$-MMP$/U$ with scaling of $A'$ terminates with a $\mathbb Q$-factorial minimal model $(X'',B'')/U$ of $(X',B')/U$. Thus the $(K_X+B)$-MMP$/U$ terminates with a minimal model $(\bar X,\bar B)/U$ of $(X,B)/U$. By \cite[Remark 2.7]{Bir12} and since $h$ is small, $(\bar X,\bar B)/U$ is a good minimal model of $(X,B)/U$ and $(X'',B'')/U$ is a $\mathbb Q$-factorial good minimal model of $(X,B)/U$. 
\end{proof}

\begin{lem}\label{lem: klt with nef and abundant boundary}
    Let $(X,B)/U$ be a klt pair and let $D$ be a nef$/U$ and abundant$/U$ $\mathbb R$-divisor on $X$. Then there exists a klt pair $(X,\Delta)$ such that
    $$K_X+\Delta\sim_{\mathbb R,U}K_X+B+D.$$
\end{lem}
\begin{proof}
    There exist a birational morphism$/U$ $h: X'\rightarrow X$ and a contraction$/U$ $g: X'\rightarrow Z$ such that $h^*D\sim_{\mathbb R,U}g^*L$ for some big$/U$ and nef$/U$ $\Rr$-divisor $L$ on $Z$. We may write $L=A_n+\frac{1}{n}E$ for some $\mathbb R$-divisor $E\geq 0$ and ample$/U$ $\mathbb R$-divisors $A_n$ for any positive integer $n$. Let $K_{X'}+B':=h^*(K_X+B)$, then there exists $n\gg 0$ such that $(X',B'+\frac{1}{n}g^*E)$ is sub-klt, so we may pick $H_n\in |A_n|_{\mathbb R/U}$ such that $(X',B'+g^*H_n+\frac{1}{n}g^*E)$ is sub-klt. We may let $\Delta:=h_*(B'+g^*H_n+\frac{1}{n}g^*E)$.
\end{proof}

\subsection{Proof of the Theorem \ref{thm: main abundance aiafs}}

The following result is the main theorem of this section which is crucial to the proof of Theorem \ref{thm: main abundance aiafs}. Its proof is similar to the proof of \cite[Theorem 5.2]{Cas+25a}. The major difference is that we use the abundance conjecture for klt pairs in dimension $d$ instead of the artificial polarization of an ample $\mathbb R$-divisor to deduce the existnece of certain good minimal models.

\begin{thm}\label{thm: log smooth existence of good minimal model}
Assume Conjecture \ref{conj: abundance} in dimension $d$. Let $\Aa:=(X,\Ff,B^{\ninv}+(1-t)B^{\inv},t)$ be a projective foliated log smooth klt algebraically integrable adjoint foliated structure of dimension $\leq d$ such that $K_{\Aa}$ is pseudo-effective. Then $\Aa$ has a $\mathbb Q$-factorial good minimal model.
\end{thm}
\begin{proof}
Let $\Aa_s:=(X,\Ff,B^{\ninv}+(1-s)B^{\inv},s)$ for any $s\in [0,1]$. By \cite[Theorem 5.1]{Cas+25a}, $K_{\Aa_s}$ is pseudo-effective for any $s\in [t,1]$. We consider the sets
$$\mathcal{P}:=\{s\in [t,1]\mid \Aa_s\text{ has a }\mathbb Q\text{-factorial good minimal model}\}$$
and
$$\mathcal{Q}:=\{s\in [t,1]\mid \Aa_s\text{ has a }\mathbb Q\text{-factorial minimal model and }K_{\Aa_s}\text{ is abundant}\}.$$
It is obvious that $\mathcal{P}\subset\mathcal{Q}$. By Theorem \ref{thm: log abundance aif}, $1\in\mathcal{Q}$. 

\medskip

\noindent\textbf{Step 1.} In this step we prove that, if $\mu\in\mathcal{Q}$ and $\mu>t$, then $\mu-\epsilon\in\mathcal{P}$ for any $0<\epsilon\ll 1$. 

Let $\Aa_{\mu}'$ be a $\mathbb Q$-factorial minimal model of $\Aa_{\mu}$. Let $\phi\colon X\dashrightarrow X'$
be the induced birational map and let $\Aa'_s:=\phi_*\Aa_s$ for any $s\in [0,1]$. Since $\phi$ is $K_{\Aa_{\mu}}$-negative, for any $0<\epsilon\ll\mu-t$ and any prime divisor $E$ on $X$, we have
$$a(E,\Aa_{\mu-\epsilon})\leq a(E,\Aa'_{\mu-\epsilon})$$
and strict inequality holds if $E$ is exceptional$/X'$.
Since $\Aa'_{\mu}$ is klt, by \cite[Proposition 3.3]{Cas+24}, $\Aa'_0$ is a klt. Since $K_{\Aa'_{\mu}}$ is nef and abundant, by Lemma \ref{lem: klt with nef and abundant boundary}, there exists a klt pair $(X',\Delta)$ such that 
$$K_{X'}+\Delta\sim_{\mathbb R}K_{\Aa'_0}+\frac{\mu-\epsilon}{\epsilon}K_{\Aa'_{\mu}}\sim_{\mathbb R,U}\frac{\mu}{\epsilon}K_{\Aa'_{\mu-\epsilon}}.$$
By Lemma \ref{lem: abundant imply gmm}, $(X',\Delta)$ has a $\mathbb Q$-factorial good minimal model $(X'',\Delta'')$ associated with birational map $\psi: X'\dashrightarrow X''$. Then $\psi_*\Aa'_{\mu-\epsilon}$ is a $\mathbb Q$-factorial good minimal model of $\Aa'_{\mu-\epsilon}$. Since $\phi$ is $K_{\Aa_{\mu-\epsilon}}$-negative, 
$$a(E,\Aa'_{\mu-\epsilon})\leq a(E,\psi_*\Aa'_{\mu-\epsilon})$$
for any prime divisor $E$ on $X'$, and strict inequality holds if $E$ is exceptional$/X''$. Therefore, $a(E,\Aa_{\mu-\epsilon})<a(E,\psi_*\Aa'_{\mu-\epsilon})$ for any prime divisor $E$ on $X$ that is exceptional$/X''$.  By Lemma \ref{lem: weak lc model only check codim 1}, $\psi\circ\phi$ is $K_{\Aa_{\mu-\epsilon}}$-negative, hence $\psi_*\Aa'_{\mu-\epsilon}$ is a $\mathbb Q$-factorial good minimal model of $\Aa_{\mu-\epsilon}$.  Therefore, $\Aa_{\mu-\epsilon}$ has a $\mathbb Q$-factorial good minimal model for any $0<\epsilon\ll\mu-t$.

\medskip

\noindent\textbf{Step 2.} In this step we prove that, if $\mu\in [t,1)$ and $\mu+\epsilon\in\mathcal{P}$ for any $0<\epsilon\ll 1-\mu$, then $\mu\in\mathcal{P}$.

Let $N_s:=N_{\sigma}(\Aa_s)$ for any $s\in [\mu,1]$. By Lemma \ref{lem: limit of nakayama-zariski decomposition}, there exists $1>s_0>\mu$ and a reduced divisor $E$ on $X$ such that $s_0-\mu\ll 1$, and $\Supp N_\mu\subset E$ and $\Supp N_s=E$ for any $s\in (\mu,s_0]$.

We let $\phi: X\dashrightarrow X'$ be a $\mathbb Q$-factorial good minimal model of $\Aa_{s_0}$ and let $\Aa'_s:=\phi_*\Aa_s$ for any $s\in [0,1]$. By Lemma \ref{lem: nz for lc divisor}, $E=\Exc(\phi)$. By Lemma \ref{lem: if contract n then movable}, $K_{\Aa_s'}$ is movable for any $s\in [\mu,s_0]$. 

Since $K_{\Aa'_{s_0}}$ is semi-ample, we may write $K_{\Aa'_{s_0}}=\sum a_iN_i$ where each $a_i>0$ and each $N_i$ is a base-point-free divisor. Let $0<\delta_0\ll s_0-\mu$ be a real number, such that $a_i\frac{s_0-\delta_0}{\delta_0}>2d$ for any $i$.

We run a $K_{\Aa'_{s_0-\delta_0}}$-MMP with scaling of an ample divisor. Since $\Aa'_{s_0}$ is klt, by \cite[Proposition 3.3]{Cas+24}, $\Aa'_0$ is a klt. Since $K_{\Aa'_{s_0}}$ is nef and abundant, by Lemma \ref{lem: klt with nef and abundant boundary}, there exists a klt pair $(X',\Delta)$ such that
$$K_{X'}+\Delta\sim_{\mathbb R}K_{\Aa'_0}+\frac{s_0-\delta_0}{\delta_0}K_{\Aa'_{s_0}}\sim_{\mathbb R}\frac{s_0}{\delta_0}K_{\Aa'_{s_0-\delta_0}}.$$
This MMP is also a $(K_{\Aa'_0}+\frac{s_0-\delta_0}{\delta_0}K_{\Aa'_{s_0}})$-MMP with scaling of an ample divisor, and is $K_{\Aa'_{s_0}}$-trivial by \cite[Lemma 4.4(3)]{BZ16}. By Lemma \ref{lem: abundant imply gmm}, this MMP terminates with a good minimal model $\Aa''_{s_0-\delta_0}$ of $\Aa'_{s_0-\delta_0}$. Let $\psi: X'\dashrightarrow X''$ be the induced birational map and let $\Aa''_s:=\psi_*\Aa'_s$ for any $s\in [0,1]$. Since $K_{\Aa'_{s_0-\delta_0}}$ is movable, by Lemma \ref{lem: nz for lc divisor}, $\psi$ does not contract any divisor, and $K_{\Aa_s''}$
is movable for any $s\in [\mu,s_0]$. Since $\psi$ is $K_{\Aa'_{s_0}}$-trivial, $K_{\Aa''_{s_0}}$ is semi-ample and $\Aa''_{s_0}$ is klt. By \cite[Proposition 3.3]{Cas+24}, $\Aa''_{0}$ is klt.

Let $B'':=(\psi\circ\phi)_*B$. Then $\Aa''_{0}=(X'',B'')$. Since
$$K_{\Aa_0''}+\frac{s}{s_0-s}K_{\Aa_{s_0}''}\sim_{\mathbb R}\frac{s_0}{s_0-s}K_{\Aa''_{s}}$$
for any $s\in [0,s_0)$ and $K_{\Aa''_{s_0}}$ is semi-ample, we may choose $0\leq H\sim_{\mathbb R}K_{\Aa_{s_0}''}$ such that
$$\left(X'',B''+\frac{s_0-\delta_0}{\delta_0}H\right)$$
is a klt pair. Since
$$K_{X''}+B''+\frac{s_0-\delta_0}{\delta_0}H\sim_{\mathbb R}\frac{s_0}{\delta_0}K_{\Aa''_{s_0-\delta_0}}$$
is nef, by Lemma \ref{lem: gmmp scaling numbers go to 0}, we may run a $(K_{X''}+B'')$-MMP with scaling of $H$ which terminates. This MMP is also a $K_{\Aa_0''}$-MMP with scaling of $K_{\Aa_{s_0}''}$.  Let $f_i: X_i\dashrightarrow X_{i+1}$ be the steps of this MMP where $X_1:=X''$ and with scaling numbers $\lambda_i$.

Let $n$ be the minimal positive integer such that $\lambda_n\leq \frac{\mu}{s_0-\mu}$. Such $n$ exists because the MMP terminates. Then there exists $\mu'>\mu$ such that the induced birational map $\tau: X''\dashrightarrow \bar X:=X_n$ is a $K_{\Aa''_s}$-MMP with scaling of $K_{\Aa''_{s_0}}$ for any $s\in [\mu,\mu')$. After decreasing $\mu'$ if necessary, we may assume that $\mu'\leq s_0-\delta_0$. Let $\bar\Aa_s:=\tau_*\Aa''_s$ for any $s\in [0,1]$. For $s\in[\mu,\mu')$, set
$$\Delta_s:=B''+\frac{s}{s_0-s}H\qquad\text{and}\qquad\bar\Delta_s:=\tau_*\Delta_s.$$
Then $(X'',\Delta_s)$ and $(\bar X,\bar\Delta_s)$ are klt, $\tau$ is a $(K_{X''}+\Delta_s)$-MMP, and, by the choice of $n$, $K_{\bar\Aa_s}$ is nef. Moreover,
$$K_{\bar X}+\bar\Delta_s\sim_{\mathbb R}\frac{s_0}{s_0-s}K_{\bar\Aa_s}.$$
By Conjecture \ref{conj: abundance}, $K_{\bar X}+\bar\Delta_s$ is abundant. Since it is already nef, Lemma \ref{lem: abundant imply gmm}, applied to the zero-step MMP of $(\bar X,\bar\Delta_s)$, shows that it is semi-ample. Hence $K_{\bar\Aa_s}$ is semi-ample for every $s\in[\mu,\mu')$. Since $K_{\Aa''_s}$ is movable for any $s\in [\mu,s_0]$, by Lemma \ref{lem: nz for lc divisor}, $\tau$ does not contract any divisor.

For any $s\in (\mu,\mu')$, let $\Aa^s_s$ be a good minimal model of $\Aa_s$ whose existence is guaranteed by our assumption. Let $\phi_s: X\dashrightarrow X^s$ be the induced birational map. By Lemma \ref{lem: nz for lc divisor}, $\Exc(\phi_s)=\Supp N_s=E$. Thus $X^s$ and $\bar X$ are isomorphic in codimension $1$. Let $p_s: W_s\rightarrow X^s$ and $q_s: W_s\rightarrow\bar X$ be a common resolution and let
$$p_s^*K_{\Aa^s_s}=q_s^*K_{\bar\Aa_s}+E_s,$$
then $E_s$ is exceptional$/X^s$ and exceptional$/\bar X$. By applying the negativity lemma twice, we have that $E_s=0$. Thus $\bar\Aa_s$ is a weak lc model of $\Aa_s$ for any $s\in (\mu,\mu')$. Let $p: W\rightarrow X$ and $q: W\rightarrow\bar X$ be a fixed common resolution. For $s\in [\mu,\mu')$, define
$$
F_s:=p^*K_{\Aa_s}-q^*K_{\bar\Aa_s}.
$$
For $s>\mu$, the weak lc property gives $F_s\geq 0$. Both terms defining $F_s$ depend affine-linearly on $s$, so
$$
F_\mu=\lim_{s\to\mu^+}F_s
$$
coefficientwise. In particular, $F_\mu\geq 0$, and
$$
p^*K_{\Aa_\mu}=q^*K_{\bar\Aa_\mu}+F_\mu.
$$

Therefore, $\bar\Aa_{\mu}$ is a semi-ample model of $\Aa_{\mu}$. By \cite[Proposition 4.2]{Cas+25a}, $\Aa_{\mu}$ has a $\mathbb Q$-factorial good minimal model. 

\medskip

\noindent\textbf{Step 3.} Since $1\in\mathcal{Q}$, by \textbf{Step 1} and \textbf{Step 2}, we have $t\in\mathcal{P}$ and we are done.
\end{proof}

\begin{thm}\label{thm: log abundance aiafs assume abundance}
Assume Conjecture \ref{conj: abundance} in dimension $d$. Let $\Aa:=(X,\Ff,B,t)$ be a projective klt algebraically integrable adjoint foliated structure of dimension $\leq d$. Then we may run a $K_{\Aa}$-MMP which scaling of an ample $\mathbb R$-divisor and any such MMP terminates with either a good minimal model or a Mori fiber space.
\end{thm}
\begin{proof}
By \cite[Theorem 1.10(1)]{Cas+24}, $X$ is potentially klt. If $K_{\Aa}$ is not pseudo-effective, then we are done by \cite[Theorem 2.1.5]{Cas+25a}. Thus we may assume that $K_{\Aa}$ is pseudo-effective. 

If $\Ff=T_X$, then we are done by Theorem \ref{thm: log abundance aif}, so we may assume that $\Ff\not=T_X$. In this case, $t<1$. Let $h: X'\rightarrow X$ be a foliated log resolution of $\Aa$ and $\Ff':=h^{-1}\Ff$. Then there exists $0<\epsilon\ll 1$ such that
$$tK_{\Ff'}+(1-t)K_{X'}+h^{-1}_*B+(1-\epsilon)(\Exc(h)^{\ninv}+(1-t)\Exc(h)^{\inv})=h^*K_{\Aa}+E$$
for some $E\geq 0$ that is exceptional$/X$ and $\Supp E=\Exc(h)$. Let $B':=h^{-1}_*B+(1-\epsilon)(\Exc(h)^{\ninv}+(1-t)\Exc(h)^{\inv})$, then $\Aa':=(X',\Ff',B',t)$ is foliated log smooth and klt, and $K_{\Aa'}$ is pseudo-effective. By Theorem \ref{thm: log smooth existence of good minimal model}, $\Aa'$ has a good minimal model. Moreover, $\Aa'$ is a foliated log smooth model of $\Aa$. By \cite[Lemma 3.8]{Cas+25b}, $\Aa$ has a bs-good minimal model. The theorem follows from \cite[Theorem 4.5]{Cas+25b}.
\end{proof}

\begin{proof}[Proof of Theorem \ref{thm: main abundance aiafs}]
    This is a special case of Theorem \ref{thm: log abundance aiafs assume abundance}.
\end{proof}

\section{Surface abundance}\label{sec: abu surface}

The goal of this section is to prove Theorem \ref{thm: abundance afs surface}. The idea is to apply the same lines of the arguments as the proof of Theorem \ref{thm: log smooth existence of good minimal model}. First, we need to prove the surface non-algebraically integrable version of \cite[Proposition 4.2, Theorem 5.1]{Cas+25a} and \cite[Proposition 3.3]{Cas+24} which were applied in the proof of Theorem \ref{thm: log smooth existence of good minimal model}.

\begin{rem}
The definition of ``foliated log smooth" for algebraically integrable foliations and for surface foliations do not coincide in literature (for algebraically integrable foliation on surfaces). To avoid confusion, in this section, and in this section only, the term ``foliated log smooth" and the related definitions, e.g. F-dlt and foliated log resolution, for rank $1$ foliations on surfaces align with the ones in \cite[Definition 4.4]{LLM23} and do not align with the definitions for algebraically integrable foliations (e.g. \cite{ACSS21}).
\end{rem}

\subsection{Surface adjoint foliated structures}\label{subsec: surface afs}

The following proposition can be seen as the surface non-algebraically integrable case of \cite[Proposition 3.3]{Cas+24}.

\begin{prop}\label{prop: klt of surface aiafs}
    Let $\Aa/U:=(X,\Ff,B^{\ninv}+(1-t)B^{\inv},t)/U$ be an lc (resp. klt) adjoint foliated structure such that $\dim X=2$ and $t<1$. Then $(X,B)$ is lc (resp. klt).
\end{prop}
\begin{proof}
We may assume that $\rk\Ff=1$, otherwise the proposition is trivial. Let $h: X'\rightarrow X$ be a foliated log resolution of $(X,\Ff,\Supp B)$ (in the sense of cf. \cite[Definition 4.4]{LLM23}). Let $\Ff':=h^{-1}\Ff$, $B':=h^{-1}_*B+\Exc(h)$, and $\Aa':=(X',\Ff',B',t)$. By \cite[Theorem 5.9]{LLM23}, we may run a $(K_{\Ff'}+B'^{\ninv})$-MMP$/X$ starting from $X'$. It terminates with a minimal model $(Y,\Ff_Y,B_Y^{\ninv})/X$ of $(X',\Ff',B'^{\ninv})/X$, with induced birational contraction $\phi:X'\rightarrow Y$; moreover, $(Y,\Ff_Y,B_Y^{\ninv})$ is F-dlt (cf. \cite[Definition 3.1]{CS21}). Let $g: Y\rightarrow X$ be the induced birational morphism and let $B_Y$ be the image of $B'$ on $Y$. Then $B_Y=g^{-1}_*B+\Exc(g)$. Moreover, since $K_{\Ff_Y}+B_Y^{\ninv}$ is nef$/X$, by the negativity lemma, for any prime divisor $D$ on $Y$ that is exceptional$/X$,
$$a(D,\Ff_Y,B_Y^{\ninv})\leq -\epsilon_{\Ff}(D).$$
Here we consider $(Y,\Ff_Y,B_Y^{\ninv})$ as a surface  foliated numerical triple (cf. \cite[Definition 3.10]{LMX24})
hence
\begin{align*}
  -t\epsilon_{\Ff}(D)-(1-t)&\leq\text{(resp. }<\text{)}a(D,X,\Ff,B^{\ninv}+(1-t)B^{\inv},t)\\
  &=ta(D,\Ff,B^{\ninv})+(1-t)a(D,X,B)\leq -t\epsilon_{\Ff}(D)+(1-t)a(D,X,B),
\end{align*}
and so $a(D,X,B)\geq -1$ (resp. $>-1$). Therefore,
$$K_Y+B_Y=K_Y+g^{-1}_*B+\Exc(g)\geq g^*(K_X+B).$$
$$\text{(resp. }K_Y+B_Y-\epsilon\Exc(g)=K_Y+g^{-1}_*B+(1-\epsilon)\Exc(g)\geq g^*(K_X+B)\text{ for some }\epsilon>0\text{)}.$$
 Since $(Y,\Ff_Y,B_Y^{\ninv})$ is dlt, $(Y,B_Y)$ is lc (resp. $(Y,B_Y-\epsilon\Exc(g))$ is klt) by \cite[Lemma 8.14]{Spi20} (note that \cite[Lemma 8.14]{Spi20} is stated for threefolds, but the same arguments, together with the proof of of \cite[Lemma 8.9]{Spi20} referred in the proof of \cite[Lemma 8.14]{Spi20} work for surfaces). Therefore, $(X,B)$ is numerically lc (resp. klt), hence lc (resp. klt).
\end{proof}

The following proposition is an analogue of \cite[Theorem 5.1]{Cas+25a}.

\begin{lem}\label{lem: afs pe surface}
Let $(X,\Ff,B,\Mm)/U$ be a $\mathbb Q$-factorial generalized foliated quadruple such that $\dim X=2$, and let $\Aa_s:=(X,\Ff,B^{\ninv}+(1-s)B^{\inv},\Mm,s)$ for any $s\in [0,1]$. Assume that $\Aa_{t}$ is lc and $K_{\Aa_{t}}$ is pseudo-effective$/U$ for some $t\in [0,1]$. Then $K_{\Aa_s}$ is pseudo-effective$/U$ for any $s\in [t,1]$.
\end{lem}
\begin{proof}
We may assume that $\rk\Ff=1$, otherwise the lemma is trivial. Assume that $K_{\Aa_1}$ is not pseudo-effective$/U$. By \cite[Theorem 1.1]{CP19}, \cite[Theorem 3.1]{LLM23}, $\Ff$ is algebraically integrable, and we get a contradiction to \cite[Theorem 5.1]{Cas+25a}.
\end{proof}

The following lemma is an analogue of \cite[Theorem 4.1]{Cas+25a}.

\begin{lem}\label{lem: extract divisors surface}
    Let $\Aa/U:=(X,\Ff,B^{\ninv}+(1-t)B^{\inv},t)/U$ be a klt adjoint foliated structure such that $\dim X=2,t<1$, and let $\mathcal{S}$ be a finite set of prime divisors over $X$ such that $a(D,\Aa)\leq 0$ for any $D\in\mathcal{S}$. Then there exists a projective birational morphism $f: Y\rightarrow X$ such that $Y$ is $\mathbb Q$-factorial, and the divisors contracted by $f$ are exactly the divisors in $\mathcal{S}$.
\end{lem}
\begin{proof}
We may assume that $\rk\Ff=1$, otherwise $a(D,X,B)\leq 0$ and we are done by \cite[Corollary 1.4.3]{BCHM10}. By Proposition \ref{prop: klt of surface aiafs}, $(X,B)$ is klt, hence $X$ is $\mathbb Q$-factorial. By induction on the number of divisors in $\mathcal{S}$, we may assume that $\mathcal{S}$ contains exactly one prime divisor $E$. If $a(E,X,B)\leq 0$, then we are done by \cite[Corollary 1.4.3]{BCHM10}, so we may assume that $a(E,X,B)>0$. Since $a(E,\Aa)\leq 0$, by linearity of discrepancies, $a(E,\Ff,B^{\ninv})<0$. For
$$
\Aa_s:=(X,\Ff,B^{\ninv}+(1-s)B^{\inv},s),
$$
we have
$$
a(E,\Aa_s)=s\,a(E,\Ff,B^{\ninv})+(1-s)a(E,X,B).
$$
Thus increasing $s$ slightly makes this discrepancy strictly negative. Since $t<1$ and the klt condition is open under a sufficiently small change of the parameter, we may choose $0<\delta_1\ll 1-t$ so that $\Aa_{t+\delta_1}$ is klt and $a(E,\Aa_{t+\delta_1})<0$. The required extraction depends only on $E$, so we relabel $t+\delta_1$ as $t$ and $\Aa_{t+\delta_1}$ as $\Aa$. Let $b:=-a(E,\Aa)$, then $b>0$.

Let $h: W\rightarrow X$ be a foliated log resolution of $\Aa$ such that $E$ is on $W$. Since $\Aa$ is klt, there exists $\delta_2>0$ such that $a(D,\Aa)>-(t\epsilon_{\Ff}(D)+(1-t))(1-2\delta_2)$
for any prime divisor $D$ on $W$ that is exceptional$/X$.

We let $B_W,F_W$ be the unique $\Rr$-divisors on $W$ such that
\begin{itemize}
\item $B_W-h^{-1}_*B:=F_W$ is exceptional$/X$,
\item $\mult_DF_W=1-\delta_2$ if $D\not=E$ and $D$ is exceptional$/X$, and 
\item $(t\epsilon_{\Ff}(E)+(1-t))\mult_EF_W=b$.
\end{itemize}
Let
$$\Aa_{W,s}:=\left(W,\Ff_W,B_W^{\ninv}+(1-s)B_W^{\inv},s\right)$$
for any $s\in [0,1]$. Then each $\Aa_{W,s}$ is foliated log smooth, and $\Aa_{W,s}$ is klt for any $s<1$. 

By \cite[Theorem 5.9]{LLM23}, we may run a $K_{\Aa_{W,1}}$-MMP$/X$ which terminates with a $\mathbb Q$-factorial minimal model $\Aa_{V,1}/X$ of $\Aa_{W,1}/X$ such that $\Aa_{V,0}$ is klt, with induced birational morphism $\phi: W\rightarrow V$ and $\Aa_{V,s}:=\phi_*\Aa_{W,s}$ for any $s\in [0,1]$. Let $g: V\rightarrow X$ be the induced birational morphism. By our construction, we have
$$K_{\Aa_{V,t}}=g^*K_{\Aa}+\sum_D((t\epsilon_{\Ff}(D)+(1-t))(1-\delta_2)+a(D,\Aa))D=: g^*K_{\Aa}+G$$
where the sum runs through all $g$-exceptional prime divisors that are not $\phi_*E$. 
By our construction of $\delta_2$, $G\geq 0$, $\Supp G=\Exc(g)$ if $E$ is contracted by $\phi$, and $\Supp G$ is the closure of $\Exc(g)\backslash\{\phi_*E\}$ if $E$ is not contracted by $\phi$. 

Let $\Pp:=\overline{K_{\Aa_{V,1}}}$. Then $\Pp$ is nef$/X$. Since $\Aa_{V,0}$ is klt, $\Bb_{V}/X:=(\Aa_{V,0},\frac{t}{1-t}\Pp)/X$ is a $\mathbb Q$-factorial generalized pair, and $K_{\Aa_{V,t}}=(1-t)K_{\Bb_V}.$ By \cite[Lemma 4.4]{BZ16} and \cite[Lemma 3.3]{Bir12}, we may run a $K_{\Bb_V}$-MMP$/X$ with scaling of an ample divisor which is also a $K_{\Aa_{V,t}}$-MMP$/X$ with scaling of an ample divisor and terminates with a model $T$. We have
$$
(1-t)K_{\Bb_V}=K_{\Aa_{V,t}}\equiv_X G.
$$
Every curve contracted by this MMP is $G$-negative over $X$. Since $V$ is a surface and $G\geq 0$, such a curve is contained in $\Supp G$; hence no divisor outside $\Supp G$ is contracted. Let $G_T$ be the pushforward of $G$ to $T$. The divisor $G_T$ is effective, exceptional$/X$, and nef$/X$, so the negativity lemma gives $G_T=0$. Hence every component of $G$ is contracted, and the divisors contracted by the induced birational map $\psi: V\dashrightarrow T$ are exactly the components of $\Supp G$.

If $E$ is not contracted by $\phi$, then since $\phi_*E$ is not a component of $\Supp G$, $E$ is also not contracted by $\psi$, hence we may let $T=:Y$ and the induced birational morphism $f: Y\rightarrow X$ satisfies our requirements. Therefore, we may assume that $E$ is contracted by $\phi$. In this case, since $X$ is $\mathbb Q$-factorial and the induced birational map $T\rightarrow X$ does not contract 
any divisor, we have $T=X$ hence $\psi: V\rightarrow X$ is a morphism. Let
$$\Bb:=\left(X,B,\frac{t}{1-t}\Pp\right),$$
then since $\Supp F_W\subset\Exc(h)$, $\Bb=\psi_*\Bb_V$, hence $\Bb$ is a klt generalized pair. We have $\Aa=\psi_*\Aa_{V,t}$.

We set $b_1:=-a(E,\Aa)$, $b_2:=-a(E,\Bb)$, $c_1:=-a(E,\Aa_{V,t})$, and $c_2:=-a(E,\Bb_V)$. Since $\psi_*K_{\Bb_V}=K_{\Bb}$, $\psi_*K_{\Aa_{V,t}}=K_{\Aa}$, and $K_{\Aa_{V,t}}=(1-t)K_{\Bb_V}$, we have
$b_1-c_1=(1-t)(b_2-c_2)$.

Let $h: Z\rightarrow V$ be a projective birational morphism such that $E$ is on $Z$. Let $\Aa_{Z,s}:=h^{-1}_*\Aa_{V,s}$ for any $s\in [0,1]$ and $\Bb_Z:=h^{-1}_*\Bb_V$.
Then we have $\mult_E(h^*K_{\Aa_{V,t}}-K_{\Aa_{Z,t}})=c_1$ and $\mult_E(h^*K_{\Bb_V}-K_{\Bb_Z})=c_2$.  Since $E$ is contracted by $\phi$, we have 
\begin{align*}
s:=&-a(E,\Aa_{V,1})
=\mult_E(h^*K_{\Aa_{V,1}}-K_{\Aa_{Z,1}})\\
<&-a(E,\Aa_{W,1})=\mult_E(B_W^{\ninv}).
\end{align*}
By linearity of discrepancies, we have
$$c_2=-a(E,\Aa_{V,0})=\frac{1}{1-t}(-a(E,\Aa_{V,t})+ta(E,\Aa_{V,1}))=\frac{1}{1-t}(c_1-ts),$$
hence, 
$c_1=(1-t)c_2+ts$.
This implies that
$$b_2=\frac{1}{1-t}(b_1-c_1)+c_2=\frac{1}{1-t}(b_1-ts).$$
By our construction of $B_W$, we have $b_1>0$ and
\[
s<\mult_E(B_W^{\ninv})\leq
\mult_E\bigl(B_W^{\ninv}+(1-t)B_W^{\inv}\bigr)=b_1.
\]
Since $0\leq t<1$, we obtain
\[
b_1-ts=(1-t)b_1+t(b_1-s)>0,
\qquad b_2=\frac{b_1-ts}{1-t}>0.
\] Since $\Bb$ is a klt generalized pair, the existence of $f$ follows from \cite[Lemma 4.6]{BZ16}.
\end{proof}

The following proposition is an analogue of \cite[Proposition 4.2]{Cas+25a}.

\begin{prop}\label{prop: eowlm implies eomm}
  Let $\Aa=(X,\Ff,B,t)/U$ be a klt adjoint foliated structure such that $\dim X=2,\rk\Ff=1$, and $t<1$. If $\Aa/U$ has a weak lc model (resp. a semi-ample model) then $\Aa/U$ has a $\mathbb Q$-factorial minimal model (resp. $\mathbb Q$-factorial good minimal model).
\end{prop}
\begin{proof}
    Let $\Aa'/U$ be a weak lc model (resp. semi-ample model) of $\Aa/U$ and let $X,X'$ be the ambient varieties of $\Aa,\Aa'$ respectively. Let 
    $$\mathcal{S}:=\{E\mid E\text{ is an exceptional}/X'\text{ prime divisor on }X, a(E,\Aa)=a(E,\Aa')\}.$$
    Then $\mathcal{S}$ is a finite set. Moreover, for any $E\in\mathcal{S}$, $a(E,\Aa')\leq 0$. By Lemma \ref{lem: extract divisors surface}, there exists a birational morphism $f: Y\rightarrow X'$ such that $Y$ is $\mathbb Q$-factorial and the divisors contracted by $f$ are exactly the divisors in $\mathcal{S}$. Let $\phi: X\dashrightarrow Y$ be the induced birational map. By Lemma \ref{lem: weak lc model only check codim 1}, $\phi_*\Aa/U$ is a $\mathbb Q$-factorial minimal model (resp. $\mathbb Q$-factorial good minimal model) of $\Aa/U$.
\end{proof}

\subsection{Proof of Theorem \ref{thm: abundance afs surface}}
\begin{proof}[Proof of Theorem \ref{thm: abundance afs surface}]
Let $\Aa_s:=(X,\Ff,s)$ for any $s\in [0,1]$. If $K_{\Aa_t}$ is not pseudo-effective, then
$$
\kappa_{\iota}(K_{\Aa_t})=\nu(K_{\Aa_t})=-\infty,
$$
and the theorem follows. Thus we may assume that $K_{\Aa_t}$ is pseudo-effective. By Lemma \ref{lem: afs pe surface}, $K_{\Aa_s}$ is pseudo-effective for any $s\in [t,1]$. In particular, $\nu:=\nu(K_{\Aa_s})$ is a constant for any $s\in (t,1)$ and $\nu\geq\nu(K_{\Aa_1})$. We consider the sets
$$\mathcal{P}:=\{s\in [t,1]\mid \Aa_s\text{ has a }\mathbb Q\text{-factorial good minimal model}\}$$  
and
$$\mathcal{Q}:=\{s\in [t,1]\mid \Aa_s\text{ has a }\mathbb Q\text{-factorial minimal model and }K_{\Aa_s}\text{ is abundant}\}.$$
It is obvious that $\mathcal{P}\subset\mathcal{Q}$. By following the same lines of the proof of Theorem \ref{thm: log smooth existence of good minimal model} verbatimly except that the following places:
\begin{itemize}
    \item We apply Proposition \ref{prop: klt of surface aiafs} instead of \cite[Proposition 3.3]{Cas+24}.
    \item We apply Proposition \ref{prop: eowlm implies eomm} instead of \cite[Proposition 4.2]{Cas+25a}.
\end{itemize}
We deduce the following:
\begin{itemize}
    \item If $\mu\in\mathcal{Q}$ and $\mu>t$, then there exists $\epsilon>0$ such that $[\mu-\epsilon,\mu)\subset\mathcal{P}$.
    \item If $\mu<1$ and $\mu+\epsilon\in\mathcal{P}$ for any $0<\epsilon\ll 1$, then $\mu\in\mathcal{P}$.
\end{itemize}
Therefore, we only need to show that $\mu\in\mathcal{Q}$ for some $\mu\in(t,1]$. By \cite[Theorem 3.1]{SS23} and \cite[Theorem 5.9]{LLM23}, there exists $0<\epsilon_0<1-t$ such that $\Aa_s$ has a $\mathbb Q$-factorial minimal model for any $s\in[1-\epsilon_0,1]$. If $K_{\Ff}$ is abundant, take $\mu=1$; if $\nu=2$, take any $\mu\in(1-\epsilon_0,1)$. We may therefore assume that $\nu\leq1$ and $K_{\Ff}$ is not abundant.

\medskip

\noindent\textbf{Step 1.} In this step we determine the invariants of $(X,\Ff)$.

Since $K_{\Ff}$ is pseudo-effective and not abundant, \cite[Lemma IV.3.1]{McQ08} gives $\nu(K_{\Ff})=1$. Apply \cite[Proposition IV.5.1 and Lemma IV.5.2]{McQ08} to the canonical model of $\Ff$, which exists by \cite[Theorem III.3.2]{McQ08}. A fibration by rational curves is excluded by pseudo-effectivity. Birational invariance of the foliated dimensions and of the irregularity of a smooth resolution then gives $\kappa(K_{\Ff})=-\infty$ and $q(X):=h^1(X,\mathcal O_X)=0$.

\medskip

\noindent\textbf{Step 2.} In this step we prove that $P\cdot K_X>0$, where $K_{\Ff}=P+N$ is the Zariski decomposition on $X$.

The divisors $P,N$ are $\mathbb Q$-divisors, $P$ is nef, $P^2=0$, $P\not\equiv0$, and $N\geq0$. Fix an ample divisor $H$. If $P\cdot K_X\leq0$, then for sufficiently large and divisible $m$, Serre duality and $(K_X-mP)\cdot H<0$ give $h^2(X,\mathcal O_X(mP))=0$. Since $q(X)=0$, Riemann--Roch gives
\[
h^0(X,\mathcal O_X(mP))\geq\chi(\mathcal O_X)-\frac m2P\cdot K_X\geq1.
\]
Since $mN$ is effective, it follows that $h^0(X,\mathcal O_X(mK_{\Ff}))\geq1$, contradicting $\kappa(K_{\Ff})=-\infty$. Therefore $P\cdot K_X>0$.

\medskip

\noindent\textbf{Step 3.} In this step we prove that $K_{\Aa_s}$ is big for $0<1-s\ll1$.

Set $L_s:=sP+(1-s)K_X$ for $s\in[0,1]$. Since $P\cdot K_X>0$ and $P\cdot H>0$, there exists $0<\delta<1$ such that for every $s\in(1-\delta,1)$,
\[
L_s^2=(1-s)\bigl(2sP\cdot K_X+(1-s)K_X^2\bigr)>0,
\qquad L_s\cdot H>0.
\]
For rational $s$ in this interval, Riemann--Roch and Serre duality imply that $L_s$ is big. For real $s$, choose rational $u<s<v$ in the same interval; the identity $L_s=\frac{v-s}{v-u}L_u+\frac{s-u}{v-u}L_v$ implies bigness. Thus $K_{\Aa_s}=L_s+sN$ is big and abundant. Choosing $\mu\in(1-\min\{\epsilon_0,\delta\},1)$ gives $\mu\in\mathcal Q\cap(t,1]$, completing the proof.
\end{proof}

\subsection*{Acknowledgements} The authors would like to thank Paolo Cascini, Junpeng Jiao, Wenhao Ou, Calum Spicer, and Roberto Svaldi for many helpful discussions. Part of this work was carried out while the authors participated in the Tianyuan Workshop: Algebraic Geometry during September 1--5, 2025; we thank the organizers Zhiyuan Li and Baohua Fu for their hospitality. The second author would like to thank Zhiyu Tian for his constant support.

\subsection*{Funding} The first author is supported by the National Key R\&D Program of China (\#2024YFA1014400). The second author is supported by start-up and matching funds (Y01286160) from the Southern University of Science and Technology (SUSTech).

\subsection*{Data availability} Data sharing does not apply to this article as no datasets were generated or analyzed during the current research.

\section{Declaration}

\subsection*{Conflict of interest} The author certifies that he has no affiliations with or involvement in any organization
or entity with any financial interest or non-financial interest in the subject matter or materials discussed in
this manuscript.


\begin{thebibliography}{99}
\bibitem[ACSS21]{ACSS21} F.~Ambro, P.~Cascini, V.~V.~Shokurov, and C.~Spicer, \textit{Positivity of the moduli part}, arXiv:2111.00423.

\bibitem[BFMT25]{BFMT25} B.~Bakker, S.~Filipazzi, M.~Mauri, and J.~Tsimerman, \textit{Baily--Borel compactifications of period images and the b-semiampleness conjecture}, arXiv:2508.19215.

\bibitem[BCHM10]{BCHM10} C.~Birkar, P.~Cascini, C.~D.~Hacon, and J.~M\textsuperscript{c}Kernan, \textit{Existence of minimal models for varieties of log general type}, J. Amer. Math. Soc. \textbf{23} (2010), no.~2, 405--468.

\bibitem[Bir12]{Bir12} C.~Birkar, \textit{Existence of log canonical flips and a special LMMP}, Publ. Math. Inst. Hautes \'{E}tudes Sci. \textbf{115} (2012), 325--368.

\bibitem[BZ16]{BZ16} C.~Birkar and D.-Q.~Zhang, \textit{Effectivity of Iitaka fibrations and pluricanonical systems of polarized pairs}, Publ. Math. IH\'{E}S \textbf{123} (2016), 283--331.

\bibitem[BM16]{BM16} F.~Bogomolov and F.~McQuillan, \textit{Rational curves on foliated varieties}, in \textit{Foliation theory in algebraic geometry}, Simons Symposia, Springer, Cham (2016), 21--51.

\bibitem[Bos01]{Bos01} J.-B.~Bost, \textit{Algebraic leaves of algebraic foliations over number fields}, Publ. Math. IH\'{E}S \textbf{93} (2001), 161--221.

\bibitem[Bru15]{Bru15} M.~Brunella, \textit{Birational geometry of foliations}, IMPA Monographs \textbf{1} (2015), Springer, Cham.

\bibitem[CP19]{CP19} F.~Campana and M.~P\u{a}un, \textit{Foliations with positive slopes and birational stability of orbifold cotangent bundles}, Publ. Math. IH\'{E}S \textbf{129} (2019), 1--49.

\bibitem[CP25]{CP25} J.~Cao and M.~P\u{a}un, \textit{Remarks on Relative Canonical Bundles and Algebraicity Criteria for Foliations in K\"ahler context}, arXiv:2502.02183.

\bibitem[CS21]{CS21} P.~Cascini and C.~Spicer, \textit{MMP for co-rank one foliations on threefolds}, Invent. Math. \textbf{225} (2021), no.~2, 603--690.

\bibitem[Cas$^+$24]{Cas+24} P.~Cascini, J.~Han, J.~Liu, F.~Meng, C.~Spicer, R.~Svaldi, and L.~Xie, \textit{Minimal model program for algebraically integrable adjoint foliated structures}, arXiv:2408.14258.

\bibitem[Cas$^+$25a]{Cas+25a} P.~Cascini, J.~Han, J.~Liu, F.~Meng, C.~Spicer, R.~Svaldi, and L.~Xie, \textit{On finite generation and boundedness of adjoint foliated structures}, arXiv:2504.10737.

\bibitem[Cas$^+$25b]{Cas+25b} P.~Cascini, J.~Liu, F.~Meng, R.~Svaldi, and L.~Xie, \textit{Variation of algebraically integrable adjoint foliated structures}, arXiv:2510.02498.

\bibitem[CS25a]{CS25a} P.~Cascini and C.~Spicer, \textit{MMP for algebraically integrable foliations}, in \textit{Higher Dimensional Algebraic Geometry: A Volume in Honor of V.~V.~Shokurov}, Cambridge University Press (2025), 69--84.

\bibitem[CS25b]{CS25b} P.~Cascini and C.~Spicer, \textit{Foliation adjunction}, Math. Ann. \textbf{391} (2025), 5695--5727.

\bibitem[CS25c]{CS25c} P.~Cascini and C.~Spicer, \textit{On the MMP for rank one foliations on threefolds}, Forum Math. Pi \textbf{13} (2025), e20, 1--38.

\bibitem[CS26]{CS26} P.~Cascini and C.~Spicer, \textit{On base point freeness for rank one foliations}, Pure Appl. Math. Q. \textbf{22} (2026), no.~2, 507--519.

\bibitem[CC25]{CC25} C.-W.~Chang and Y.-A.~Chen, \textit{Boundedness of toric foliations}, arXiv:2502.11080.

\bibitem[CM24]{CM24} P.~Chaudhuri and R.~Mascharak, \textit{Log canonical minimal model program for corank one foliations on threefolds}, arXiv:2410.05178.

\bibitem[CHLX23]{CHLX23} G.~Chen, J.~Han, J.~Liu, and L.~Xie, \textit{Minimal model program for algebraically integrable foliations and generalized pairs}, arXiv:2309.15823.

\bibitem[Cho08]{Cho08} S.~R.~Choi, \textit{The geography of log models and its applications}, Ph.D. thesis, Johns Hopkins University, 2008.

\bibitem[Dar78a]{Dar78a} G.~Darboux, \textit{M\'emoire sur les \'equations diff\'erentielles alg\'ebriques du premier ordre et du premier degr\'e}, Bulletin des Sciences Math\'ematiques et Astronomiques \textbf{2} (1878), no.~1, 151--200.

\bibitem[Dar78b]{Dar78b} G.~Darboux, \textit{M\'emoire sur la th\'eorie des coordonn\'ees curvilignes, et des syst\`emes orthogonaux}, Ann. Sci. \'Ecole Norm. Sup. \textbf{7} (1878), no.~2, 101--150, 227--260, 275--348.

\bibitem[Dru18]{Dru18} S.~Druel, \textit{A decomposition theorem for singular spaces with trivial canonical class of dimension at most five}, Invent. Math. \textbf{211} (2018), no.~1, 245--296.

\bibitem[Dru21]{Dru21} S.~Druel, \textit{Codimension 1 foliations with numerically trivial canonical class on singular spaces}, Duke Math. J. \textbf{170} (2021), no.~1, 95--203.

\bibitem[DO22]{DO22} S.~Druel and W.~Ou, \textit{Codimension one foliations with numerically trivial canonical class on singular spaces II}, Int. Math. Res. Not. (2022), no.~20, 15574--15611.

\bibitem[FKX17]{FKX17} T.~de Fernex, J.~Koll\'ar, and C.~Xu, \textit{The dual complex of singularities}, in \textit{Higher dimensional algebraic geometry: in honor of Professor Yujiro Kawamata's sixtieth birthday}, Adv. Stud. Pure Math. \textbf{74} (2017), Math. Soc. Japan, Tokyo, 103--129.

\bibitem[FS25]{FS25} S.~Filipazzi and C.~Spicer, \textit{On semi-ampleness of the moduli part}, Moduli \textbf{2} (2025), e12, 1--27.

\bibitem[Ghy00]{Ghy00} E.~Ghys, \textit{\`A propos d'un th\'eor\`eme de J.-P.~Jouanolou concernant les feuilles ferm\'ees des feuilletages holomorphes}, Rend. Circ. Mat. Palermo \textbf{49} (2000), 175--180.

\bibitem[GL13]{GL13} Y.~Gongyo and B.~Lehmann, \textit{Reduction maps and minimal model theory}, Compos. Math. \textbf{149} (2013), no.~2, 295--308.

\bibitem[HL23]{HL23} C.~D.~Hacon and J.~Liu, \textit{Existence of flips for generalized lc pairs}, Camb. J. Math. \textbf{11} (2023), no.~4, 795--828.

\bibitem[HJLL24]{HJLL24} J.~Han, J.~Jiao, M.~Li, and J.~Liu, \textit{Volume of algebraically integrable foliations and locally stable families}, arXiv:2406.16604.

\bibitem[HLX25]{HLX25} J.~Han, J.~Liu, and Q.~Xue, \textit{On the equivalence between the effective adjunction conjectures of Prokhorov--Shokurov and of Li}, Algebra Number Theory \textbf{19} (2025), no.~11, 2261--2279.

\bibitem[HLS26]{HLS26} J.~Han, J.~Liu, and V.~V.~Shokurov, \textit{ACC for minimal log discrepancies of exceptional singularities}, Peking Math. J. \textbf{9} (2026), no.~2, 225--257.

\bibitem[HH20]{HH20} K.~Hashizume and Z.~Hu, \textit{On minimal model theory for log abundant lc pairs}, J. Reine Angew. Math. \textbf{767} (2020), 109--159.

\bibitem[HP19]{HP19} A.~H\"oring and T.~Peternell, \textit{Algebraic integrability of foliations with numerically trivial canonical bundle}, Invent. Math. \textbf{216} (2019), 395--419.

\bibitem[Hu20]{Hu20} Z.~Hu, \textit{Log abundance of the moduli b-divisors of lc-trivial fibrations}, arXiv:2003.14379.

\bibitem[Jou78]{Jou78} J.-P.~Jouanolou, \textit{Hypersurfaces solutions d'une \'equation de Pfaff analytique}, Math. Ann. \textbf{232} (1978), no.~3, 239--245.

\bibitem[KMM87]{KMM87} Y.~Kawamata, K.~Matsuda, and K.~Matsuki, \textit{Introduction to the minimal model problem}, in \textit{Algebraic geometry, Sendai, 1985}, Adv. Stud. Pure Math. \textbf{10}, North-Holland, Amsterdam (1987), 283--360.

\bibitem[Kaw92]{Kaw92} Y.~Kawamata, \textit{Abundance theorem for minimal threefolds}, Invent. Math. \textbf{108} (1992), 229--246.

\bibitem[KMM94]{KMM94} S.~Keel, K.~Matsuki, and J.~M\textsuperscript{c}Kernan, \textit{Log abundance theorem for threefolds}, Duke Math. J. \textbf{75} (1994), 99--119.

\bibitem[Kol13]{Kol13} J.~Koll\'ar, \textit{Singularities of the minimal model program}, Cambridge Tracts in Math. \textbf{200} (2013), Cambridge Univ. Press.

\bibitem[Kol23]{Kol23} J.~Koll\'ar, \textit{Families of varieties of general type}, Cambridge Tracts in Math. \textbf{231} (2023), Cambridge Univ. Press.

\bibitem[KM98]{KM98} J.~Koll\'ar and S.~Mori, \textit{Birational geometry of algebraic varieties}, Cambridge Tracts in Math. \textbf{134} (1998), Cambridge Univ. Press.

\bibitem[KP17]{KP17} S.~J.~Kov\'acs and Z.~Patakfalvi, \textit{Projectivity of the moduli space of stable log-varieties and subadditivity of log-Kodaira dimension}, J. Amer. Math. Soc. \textbf{30} (2017), no.~4, 959--1021.

\bibitem[Lai11]{Lai11} C.-J.~Lai, \textit{Varieties fibered by good minimal models}, Math. Ann. \textbf{350} (2011), no.~3, 533--547.

\bibitem[Li25]{Li25} M.~Li, \textit{MMP for generalized foliated threefolds of rank one}, Manuscripta Math. \textbf{176} (2025), Paper No.~78.

\bibitem[LLM23]{LLM23} J.~Liu, Y.~Luo, and F.~Meng, \textit{On global ACC for foliated threefolds}, Trans. Amer. Math. Soc. \textbf{376} (2023), no.~12, 8939--8972.

\bibitem[LMX24]{LMX24} J.~Liu, F.~Meng, and L.~Xie, \textit{Complements, index theorem, and minimal log discrepancies of foliated surface singularities}, Eur. J. Math. \textbf{10} (2024), no.~1, Paper No.~6, 29 pp.

\bibitem[LMX25]{LMX25} J.~Liu, F.~Meng, and L.~Xie, \textit{Minimal model program for algebraically integrable foliations on klt varieties}, Compos. Math. \textbf{161} (2025), no.~12, 3213--3276.

\bibitem[LX25a]{LX25a} J.~Liu and L.~Xie, \textit{Relative Nakayama--Zariski decomposition and minimal models of generalized pairs}, Peking Math. J. \textbf{8} (2025), no.~2, 299--349.

\bibitem[LX25b]{LX25b} J.~Liu and Z.~Xu, \textit{Non-vanishing implies numerical dimension one abundance}, arXiv:2505.05250.

\bibitem[LPT18]{LPT18} F.~Loray, J.~Pereira, and F.~Touzet, \textit{Singular foliations with trivial canonical class}, Invent. Math. \textbf{213} (2018), 1327--1380.

\bibitem[LW24]{LW24} J.~Lu and X.-H.~Wu, \textit{On the 1-adjoint canonical divisor of a foliation}, Manuscripta Math. \textbf{175} (2024), no.~3--4, 739--752.

\bibitem[LWX26]{LWX26} J.~Lu, X.-H.~Wu, and S.~Xu, \textit{Canonical models of adjoint foliated structures on surfaces}, J. Lond. Math. Soc. \textbf{114} (2026), no.~1, Paper No.~e70622.

\bibitem[LT24]{LT24} X.~L\"u and S.~Tan, \textit{The Poincar\'e Problem for a foliated surface}, arXiv:2404.16293.

\bibitem[McQ08]{McQ08} M.~McQuillan, \textit{Canonical models of foliations}, Pure Appl. Math. Q. \textbf{4} (2008), no.~3, Special Issue: In honor of Fedor Bogomolov, Part 2, 877--1012.

\bibitem[McQ24]{McQ24} M.~McQuillan, \textit{Semi-Stable Reduction of Foliations}, in Y.~Tschinkel (ed.), \textit{Arithmetic and Algebraic Geometry: A Mathematical Tribute to Yuri Manin}, Simons Symposia, Springer, Cham (2024), 239--351.

\bibitem[Miy87]{Miy87} Y.~Miyaoka, \textit{Deformations of a morphism along a foliation and applications}, in \textit{Algebraic Geometry, Bowdoin, 1985}, Proc. Sympos. Pure Math. \textbf{46}, Part~1, Amer. Math. Soc., Providence, RI (1987), 245--268.

\bibitem[Nak04]{Nak04} N.~Nakayama, \textit{Zariski-decomposition and abundance}, MSJ Mem. \textbf{14} (2004), Math. Soc. Japan, Tokyo.

\bibitem[Ou25a]{Ou25a} W.~Ou, \textit{A characterization of uniruled compact K\"ahler manifolds}, arXiv:2501.18088.

\bibitem[Ou25b]{Ou25b} W.~Ou, \textit{Foliations whose first Chern class is nef}, J. Eur. Math. Soc. \textbf{27} (2025), 2953--2984.

\bibitem[Pai92]{Pai92} P.~Painlev\'e, \textit{M\'emoire sur les \'equations diff\'erentielles du premier ordre (suite)}, Ann. Sci. \'Ecole Norm. Sup. (3) \textbf{9} (1892), 101--144.

\bibitem[Pas24]{Pas24} A.~Passantino, \textit{Numerical conditions for the boundedness of foliated surfaces}, arXiv:2412.05986.

\bibitem[PX17]{PX17} Z.~Patakfalvi and C.~Xu, \textit{Ampleness of the CM line bundle on the moduli space of canonically polarized varieties}, Algebr. Geom. \textbf{4} (2017), no.~1, 29--39.

\bibitem[PS19]{PS19} J.~V.~Pereira and R.~Svaldi, \textit{Effective algebraic integration in bounded genus}, Algebr. Geom. \textbf{6} (2019), no.~4, 454--485.

\bibitem[Poi85]{Poi85} H.~Poincar\'e, \textit{Sur les \'equations lin\'eaires aux diff\'erentielles ordinaires et aux diff\'erences finies}, Amer. J. Math. \textbf{7} (1885), no.~3, 203--258.

\bibitem[Poi91]{Poi91} H.~Poincar\'e, \textit{Sur l'int\'egration alg\'ebrique des \'equations diff\'erentielles du premier ordre et du premier degr\'e}, Rendiconti del Circolo Matematico di Palermo (1884--1940) \textbf{5} (1891), no.~1, 161--191.

\bibitem[Sho92]{Sho92} V.~V.~Shokurov, \textit{Threefold log flips}, Izv. Ross. Akad. Nauk Ser. Mat. \textbf{56} (1992), no.~1, 105--203.

\bibitem[Spi20]{Spi20} C.~Spicer, \textit{Higher dimensional foliated Mori theory}, Compos. Math. \textbf{156} (2020), no.~1, 1--38.

\bibitem[SS22]{SS22} C.~Spicer and R.~Svaldi, \textit{Local and global applications of the Minimal Model Program for co-rank 1 foliations on threefolds}, J. Eur. Math. Soc. \textbf{24} (2022), no.~11, 3969--4025.

\bibitem[SS23]{SS23} C.~Spicer and R.~Svaldi, \textit{Effective generation for foliated surfaces: Results and applications}, J. Reine Angew. Math. \textbf{795} (2023), 45--84.

\bibitem[Tan96]{Tan96} S.~Tan, \textit{On the invariants of base changes of pencils of curves, II}, Math. Z. \textbf{222} (1996), no.~4, 655--676.

\bibitem[Tan10]{Tan10} S.~Tan, \textit{Chern numbers of a singular fiber, modular invariants and isotrivial families of curves}, Acta Math. Vietnam. \textbf{35} (2010), no.~1, 159--172.
\end{thebibliography}
\end{document}